# Rigidity of infinite disk patterns

By  Zheng-Xu He*

## Abstract

Let $P$ be a locally finite disk pattern on the complex plane $\mathbb{C}$ whose combinatorics are described by the one-skeleton $G$ of a triangulation of the open topological disk and whose dihedral angles are equal to a function $\Theta : E \to [0, \pi/2]$ on the set of edges. Let $P^*$ be a combinatorially equivalent disk pattern on the plane with the same dihedral angle function. We show that $P$ and $P^*$ differ only by a euclidean similarity.

In particular, when the dihedral angle function $\Theta$ is identically zero, this yields the rigidity theorems of B. Rodin and D. Sullivan, and of O. Schramm, whose arguments rely essentially on the pairwise disjointness of the interiors of the disks. The approach here is analytical, and uses the maximum principle, the concept of vertex extremal length, and the recurrency of a family of electrical networks obtained by placing resistors on the edges in the contact graph of the pattern.

A similar rigidity property holds for locally finite disk patterns in the hyperbolic plane, where the proof follows by a simple use of the maximum principle. Also, we have a uniformization result for disk patterns.

In a future paper, the techniques of this paper will be extended to the case when $0 \le \Theta < \pi$. In particular, we will show a rigidity property for a class of infinite convex polyhedra in the 3-dimensional hyperbolic space.

## 1. Introduction

A *disk pattern* on the Riemann sphere $\hat{\mathbb{C}}$ will be defined to be a collection of closed disks in $\hat{\mathbb{C}}$ in which no disk has its boundary contained in the union of two other disks and no disk is the Hausdorff limit of a sequence of distinct disks. The *contact graph* $G = G_P$ of such a pattern $P$ is the graph whose vertices correspond to the disks of the pattern, and an edge appears in $G$ if

*This research is partially supported by an NSF grant. The paper was circulated under the title, "Rigidity of disk packings with overlappings, I."

the corresponding disks intersect each other. The *dihedral angle* of a pair of intersecting disks $D_1$ and $D_2$ is defined to be the angle in $[0, \pi)$ between the clockwise tangent of $\partial D_1$ and the counterclockwise tangent of $\partial D_2$ at a point of $\partial D_1 \cap \partial D_2$. Let $E$ be the set of edges in the graph $G$. For any $[v_1, v_2] \in E$, let $\Theta_P([v_1, v_2])$ be the dihedral angle of the disks $P(v_1)$ and $P(v_2)$. Then $\Theta_P : E \to [0, \pi)$ is called the *dihedral angle function* of $P$.

Let $G$ be a given graph and let $\Theta : E \to [0, \pi)$ be a function defined on the set of edges of $G$. Let us consider the following problem: Does there exist a disk pattern $P$ whose contact graph is (combinatorially equivalent to) $G$ and whose dihedral angle function is $\Theta$? And if it does, to what extent is the disk pattern unique? This problem is well posed under the condition that $0 \leq \Theta \leq \pi/2$. In this case, by Thurston's interpretation of Andreev's theorem, the existence problem has a complete answer when $G$ is finite and the uniqueness is understood when $G$ is the one-skeleton of a triangulation of the 2-sphere (see e.g. [26] or [16] for precise statements). In this paper, we will study infinite disk patterns under the condition $0 \leq \Theta \leq \pi/2$. The general case $0 \leq \Theta < \pi$ is technically more involved, and will be considered in [10] where we will also generalize a characterization theorem due to I. Rivin and C. Hodgson [18] for convex polyhedra in hyperbolic 3-space.

For a vertex $v$ in the contact graph $G$ of a disk pattern $P$, we will denote $P(v)$ to be the disk in $P$ corresponding to the vertex. If $P$ is a disk pattern in the plane $\mathbb{C}$ whose dihedral angle function is bounded by $\pi/2$, then the euclidean center $A(v)$ of a disk $P(v)$ lies outside of any other disk; for an edge $[v_1, v_2]$ in $G$, the straight arc $A(v_1)A(v_2)$ is disjoint from the interior of any disk $P(v)$ with $v \neq v_1$ and $v \neq v_2$. The natural map which maps an edge $[v_1, v_2]$ homeomorphically onto the arc $A(v_1)A(v_2)$ is an immersion of the graph $G$ into $\mathbb{C}$. The only possible double points of this immersion can be described as follows: The arc $A(v_0)A(v_2)$ intersects $A(v_1)A(v_3)$ (as in Figure 1.1) if and only if: $\Theta_P([v_{i-1}, v_i]) = \pi/2$, $1 \leq i \leq 4$ (where $v_4 = v_0$), and either $\Theta_P([v_0, v_2])$ or $\Theta_P([v_1, v_3])$ is equal to 0. If this is the case, we must have both $\Theta_P([v_0, v_2]) = \Theta_P([v_1, v_3]) = 0$ and we will call the edges $[v_0, v_2]$ and $[v_1, v_3]$ *reducible*. The *reduced graph* of $P$, consisting of the same vertex set as $G$ and the irreducible edges, is embedded in $\mathbb{C}$.

Similar observations can be made for a disk pattern on the Riemann sphere, under the additional assumption that each disk in the pattern is smaller than a hemisphere (see [26]).

Let $G$ be a graph, and let $\Theta : E \to [0, \pi/2]$ be a function on the set of edges. Whenever a simple loop $v_1$, $v_2$, $v_3$, $v_4 = v_0$ in $G$ has the property that $\Theta([v_{i-1}, v_i]) = \pi/2$, $i = 1, \ldots, 4$, and $\Theta([v_0, v_2]) = 0$ as in Figure 1.1(b), then we will add the other reducible edge $[v_1, v_3]$, in case it is not in $G$, to the graph. Let $\tilde{G}$ denote the graph thus obtained, and define a function $\tilde{\Theta}$ by letting $\tilde{\Theta}(e) = \Theta(e)$ if $e$ is an edge in $G$, and $\tilde{\Theta}(e) = 0$ if $e$ is in $\tilde{G} - G$. A disk pattern $P$ is said to *realize the data* $(G, \Theta)$ if its contact graph is combinatorially isomorphic to $\tilde{G}$ and the corresponding dihedral angle function is equal to $\tilde{\Theta}$.

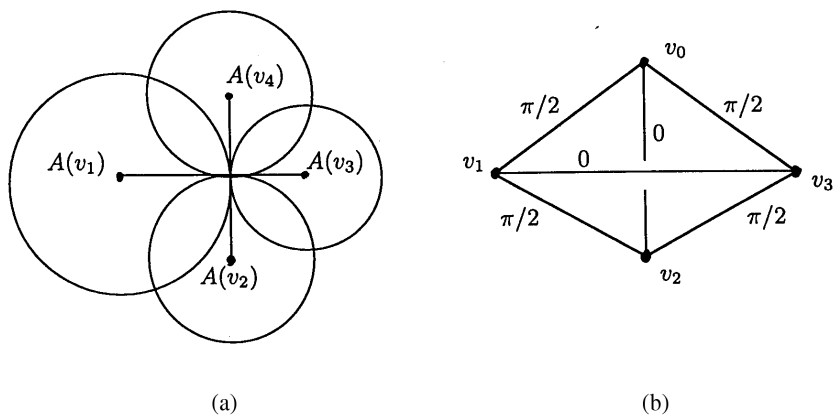

(a)                  (b)

**Figure 1.1.** (a) The configuration of the disks $P(v_i)$, $1 \leq i \leq 4$;
(b) The corresponding contact graph with dihedral angles marked on the edges.

A disk pattern $P$ is called *locally finite* in a domain $\Omega \subseteq \hat{\mathbb{C}}$, if all disks of $P$ are contained in $\Omega$ and every compact subset of $\Omega$ intersects only finitely many disks of $P$. A graph is called a *disk triangulation graph* if it is equal to the one-skeleton of a triangulation $T$ of the open topological disk. The main result of this paper is the following rigidity theorem.

RIGIDITY THEOREM 1.1. *Let $G$ be a disk triangulation graph, and let $\Theta : E \to [0, \pi/2]$ be a function defined on the set of edges. Let $P$ and $P^*$ be disk patterns in $\mathbb{C}$ which realize $(G, \Theta)$. Assume that $P$ is locally finite in the plane. Then there is a euclidean similarity $f : \mathbb{C} \to \mathbb{C}$ such that $P^* = f(P)$.*

See Figure 1.2 for several examples of locally finite disk patterns in $\mathbb{C}$. By Theorem 1.1, the geometry (i.e. the similarity class) of any of these patterns is uniquely determined by its contact graph and dihedral angle function.

When $\Theta = 0$ on all edges, then $P$ is a disk packing. Rigidity property in this case has been proved by B. Rodin and D. Sullivan [21] in the bounded valence case, using Schottky groups (see [9, Remark (p. 407)] for a simplified proof), and by O. Schramm [22] in the unbounded valence case, using an argument which is topological in nature. A more direct proof of the rigidity result of O. Schramm using a similar idea can be found in [11]. A different proof based on the Schottky groups and a "generalized Gröstch argument" is given in [12]. All these proofs rely essentially on the property that pairs of disks have disjoint interiors, and therefore cannot be generalized.

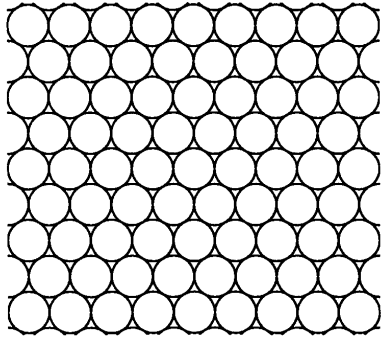

(a) A regular hexagonal packing.

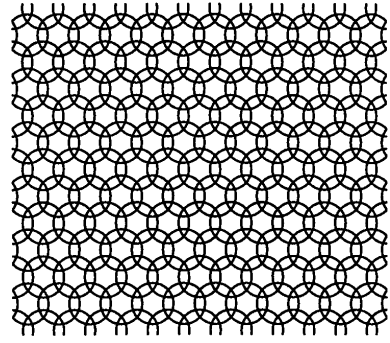

(b) A hexagonal pattern.

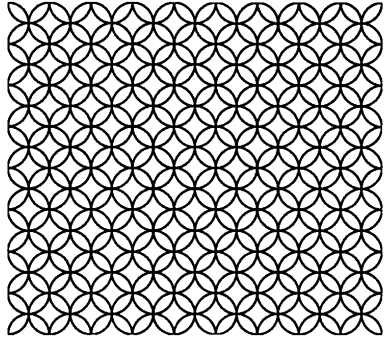

(c) A square grid pattern.

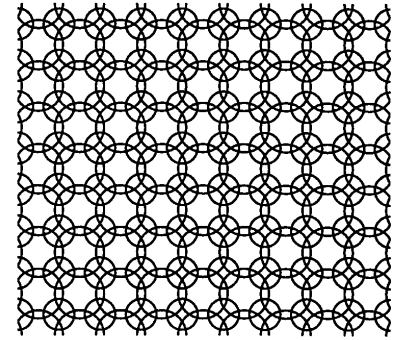

(d) A square grid with deformed angles.

**Figure 1.2.** The similarity class of any of these patterns is uniquely determined by its contact graph and dihedral angle function.

The method in the present paper is analytical, and can be roughly described as follows. Let $P$ and $P^*$ be the patterns in Rigidity Theorem 1.1. The first step in our approach is to show that there are some lower and upper bounds on the ratio of the radii of pairs of corresponding disks in $P^*$ and $P$. The bounds will depend on the patterns but not on the individual disks. This uses the idea of vertex extremal length (cf. [6], [7], [23]) as well as the maximum principle. Next we will construct a one-parameter family of "immersed disk patterns" $P_t$, $0 \le t \le 1$, with $P_0 = P$, $P_1 = P^*$. For each vertex $v$ in the graph, let $l(v, t) = \log \rho(P_t(v)) - \log \rho(P_0(v))$, where $\rho(\cdot)$ denotes the euclidean radius. Using the maximum principle again, it will be shown that the rate of change $h(v, t) = d\, l(v, t)/dt = \dfrac{d\, \rho(P_t(v))/dt}{\rho(P_t(v))}$ is bounded by a constant independent of $v$. Moreover, for each fixed $t$, $h(v, t)$ is harmonic in the electrical network obtained from the graph $G$ by placing resistors of certain conductances on its edges. In fact, the conductance of an edge $[v_0, v_1]$ is defined to be the partial derivative of the "curvature" at $v_0$ with respect to $\log(\rho(P_t(v_1)))$. Then we prove that the network is recurrent in the sense that the resistance from a finite subset of vertices to the infinity is infinite. As a bounded harmonic function, $h(v, t)$ is therefore constant in the graph. It follows that $l(v, t)$ and hence the ratio $\rho(P_t(v))/\rho(P_0(v))$ are independent of $v$. We deduce that the patterns $P_t$, in particular $P_1 = P^*$, are just images of $P_0 = P$ under euclidean similarities.

In case of disk patterns in the Poincaré disk $U = \{z \in \mathbb{C}; |z| < 1\}$, the following can be easily proved using the maximum principle in the hyperbolic plane.

RIGIDITY THEOREM 1.2. *Let $G$ be a disk triangulation graph, and let $\Theta : E \to [0, \pi/2]$ be a function defined on the set of edges. Let $P$ and $P^*$ be disk patterns in $U$ which realize $(G, \Theta)$. Assume that both $P$ and $P^*$ are locally finite in $U$. Then there is a hyperbolic isometry $f : U \to U$ such that $P^* = f(P)$.*

Let us conclude this section by considering the Koebe uniformization of disk patterns. Let $G$ be a connected, planar graph. Assume that each edge has distinct ends, and that through a pair of vertices there is at most one edge. A subset contained in $G$ is said to *separate* the vertices of $G$, if there is a pair of vertices in the complement of the subset so that any path joining the vertices passes through the subset. Let $\Theta : E \to [0, \pi/2]$ be a function on the set of edges. It is easy to show that if there is a disk pattern in $\hat{\mathbb{C}}$ which realizes the data $(G, \Theta)$, then the following conditions must hold (see [26]):

(C1) If a simple loop in $G$ formed by three edges $e_1, e_2, e_3$ separates the vertices of $G$, then $\sum_{i=1}^{3} \Theta(e_i) < \pi$;

(C2) If $v_1, v_2, v_3, v_4 = v_0$ are distinct vertices in $G$ and if $[v_{i-1}, v_i] \in E$ and $\Theta([v_{i-1}, v_i]) = \pi/2$, $i = 1, \ldots, 4$, then either $[v_0, v_2]$ or $[v_1, v_3]$ is an edge in $G$.

We note that conditions (C1), (C2) and $\Theta \leq \pi/2$ imply that if a simple loop formed by the edges $e_1, e_2, \ldots, e_k$ separates the vertices of $G$, then $\sum_{i=1}^{k} (\pi - \Theta(e_i)) > 2\pi$.

On the other hand, using Thurston's interpretation of Andreev's theorem and a compactness argument, one may show that conditions (C1) and (C2) are also sufficient for the existence of a disk pattern which realizes the data, provided that the graph has at least five vertices (see [26] and §7). In the special case of disk triangulation graphs, we have the following theorem, proved along the same lines as Theorem 7.1 of [13] (also compare [24]).

UNIFORMIZATION THEOREM 1.3. *Let $G$ be a disk triangulation graph, and let $\Theta : E \to [0, \pi/2]$ be a function defined on the set of edges. Assume that conditions* (C1) *and* (C2) *hold.*

(i) *If $G$ is* VEL-*parabolic, then there is a locally finite disk pattern in $\mathbb{C}$ which realizes the data $(G, \Theta)$.*

(ii) *If $G$ is* VEL-*hyperbolic, then there is a locally finite disk pattern in $U$ which realizes $(G, \Theta)$.*

When $\Theta = 0$, the above theorem implies Corollary 0.5 of [11], which in turn, generalizes the result of [4] in the bounded valence case. Here we also note that I. Rivin has a related uniformization theorem of finite ideal polyhedra in hyperbolic 3-space (see [17]). As we noted earlier, the characterization of compact polyhedra in terms of their polars was beautifully accomplished in the joint work [18].

The paper is organized as follows. We will start in Section 2 by presenting the Maximum Principle for disk patterns. In particular, Theorem 1.2 will be obtained as a corollary. In Section 3, we will discuss immersed disk patterns, and we will show that one may prescribe the radii of the boundary disks of some finite pattern. In Section 4, assuming that $\rho(P^*(v))/\rho(P(v))$ is uniformly bounded, we construct the variation $P_t$ from $P$ to $P^*$. We then show that for each $t$, the rate of change $h(v, t)$ of $\log \rho(P_t(v))$ is a bounded harmonic function in the electrical network obtained from the graph $G$ by placing resistors of certain conductances on its edges. Then in Section 5, we show that the network is recurrent, and as a consequence, $h(v, t)$ is independent of $v$. The uniform boundedness of $\rho(P^*(v))/\rho(P(v))$ will be proved in Section 6. Finally, in Section 7, we will study the uniformization of disk patterns.

*Notational conventions.* Throught the paper, $G$ will denote a graph, $E = E(G)$ the set of edges in $G$, and $V = V(G)$ the set of vertices. The symbol $\Theta$ will denote a function from $E$ to $[0, \pi)$. If $P$ is a disk pattern which realizes $(G, \Theta)$, then for any subset of vertices $W \subseteq V$, we denote $P(W) = \bigcup_{v \in W} P(v)$. For a circle $c$, we denote by $V_c(P)$ the set of vertices $v$ for which $P(v) \cap c \neq \emptyset$. For any $r > 0$, we denote by $c(r)$ the circle of radius $r$ centered at 0, and $D(r)$ the closed disk bounded by $c(r)$. For a disk $D$, we will denote $\rho(D)$ to be its euclidean radius.

We are very grateful to Professor Burt Rodin for introducing us to some basic techniques for the rigidity problem. We are very much indebted to Professors Michael Freedman, and Dov Aharonov, Peter Doyle, Oded Schramm and Lihe Wang for helpful discussions as well as overall influence. We are also very grateful to the referee for carefully reading the original manuscript and providing many helpful suggestions.

## 2. The maximum principle

Various forms of the maximum principle, some weaker, some stronger, have been suggested and used by some authors. Here we will present it in a way suitable to our applications. The proofs, included here for completeness, are elementary and mostly well known.

Let $G = (V, E)$ and $\Theta : E \to [0, \pi/2]$ be given. A vertex $v_0$ in $G$ is called *interior* if there is a closed chain of neighboring vertices $v_1, v_2, \ldots, v_l$, where either $l \geq 4$, or $l = 3$ and each edge $[v_0, v_k]$ is irreducible, $1 \leq k \leq 3$. Otherwise, it is called a *boundary* vertex. Note that in the graph in Figure 1.1(b), with the values of $\Theta$ marked on the edges, $v_0$ is the neighbor of the chain of vertices $v_1, v_2, v_3$, but is not an interior vertex by our definition. If $G$ is the one-skeleton of a triangulation of a planar surface possibly with boundary so that each boundary component has at least four vertices, and if (C1) and (C2) hold, then a vertex is interior precisely when it lies in the interior of the surface.

LEMMA 2.1 (maximum principle). *Let $G$ be a finite graph, and assume $0 \leq \Theta \leq \pi/2$. Let $P$ and $P^*$ be disk patterns in $\mathbb{C}$ which realize $(G, \Theta)$. Then the maximum (or minimum) of $\rho(P^*(v))/\rho(P(v))$ is attained at a boundary vertex.*

*Proof.* Let $v_0$ be an interior vertex. Let $v_1, v_2, \ldots, v_l, v_{l+1} = v_1$ be the chain of neighboring vertices. Let $Q$ be a disk pattern which realizes the data $(G, \Theta)$. Denote by $A_j$ be the center of the disk $Q(v_j)$, $0 \leq j \leq l + 1$. Let $\beta_k(Q)$ be the angle $\angle A_k A_0 A_{k+1}$. Since $v_0$ is interior, we have $\sum_{k=1}^{l} \beta_k(Q) = 2\pi$. On the other hand, $\beta_k(Q)$ is a nondecreasing function of $\rho(Q(v_k))$ (and of $\rho(Q(v_{k+1}))$) (see e.g. [26], or [16, Lemma 3 (p. 111)]). The sum $\beta_{k-1}(Q) + \beta_k(Q)$, and hence $\sum_{j=1}^{l} \beta_j(Q)$, are strictly increasing in $\rho(Q(v_k))$ if the edge $[v_0, v_k]$ is irreducible. The maximum principle follows immediately. $\square$

*Remark.* Assume that the set of interior vertices of $G$ in Lemma 2.1 is connected. If the maximum (or minimum) of $\rho(P^*(v))/\rho(P(v))$ is attained at an interior vertex, then it is constant in the union of the set of interior vertices and the set of boundary vertices which share an irreducible edge with an interior vertex.

The maximum principle suggests that the logarithm of the ratio $\rho(P^*(v))/\rho(P(v))$, considered as a function on the vertex set of $G$, behaves as a harmonic function in the graph, an idea which we will pursue in Section 4.

A version of the maximum principle also holds for patterns in the hyperbolic disk $U$. For a closed disk $D$ in $U$, denote by $\rho_{\text{hyp}}(D)$ its hyperbolic radius. Let $D_0$ be a fixed closed disk in $U$ centered at 0, and let $D_1$ be a variable closed disk which intersects $D_0$ with a fixed dihedral angle between 0 and $\pi/2$. Then $\rho_{\mathrm{hyp}}(D_1)$ is a strictly increasing function of $\rho(D_1)$, and for any real $\gamma > 1$, we have: $\rho_{\mathrm{hyp}}(\gamma D_1)/\rho_{\mathrm{hyp}}(D_1) > \rho_{\mathrm{hyp}}(\gamma D_0)/\rho_{\mathrm{hyp}}(D_0)$, where $\gamma D = \{\gamma z; z \in D\}$. Combining these facts with Lemma 2.1, we deduce:

LEMMA 2.2 (maximum principle in the hyperbolic plane). *Let $G$, $\Theta$, $P$ and $P^*$ be as in Lemma 2.1. Assume that the disks of $P$ and $P^*$ are contained in $U$. Then*:

(a) *The maximum of $\rho_{\mathrm{hyp}}(P^*(v))/\rho_{\mathrm{hyp}}(P(v))$, if $> 1$, is never attained at an interior vertex; and*

(b) *In particular, if the inequality $\rho_{\mathrm{hyp}}(P^*(v)) \leq \rho_{\mathrm{hyp}}(P(v))$ is true for each boundary vertex, then it holds for all vertices of $G$.*

The second part of the lemma is the analog of the Schwarz-Pick lemma. See [5] and [11] for versions of Schwarz-Pick lemma for disk packings (also compare [19], [20]). In [10], this will be generalized to the case $0 \leq \Theta < \pi$.

To extend the lemma, we make the following definition of $\rho_{\mathrm{hyp}}(D)$ for a closed disk $D$ in $\hat{\mathbb{C}}$ which intersects $U$. If $D$ in contained in $U$, we define $\rho_{\mathrm{hyp}}(D)$ to be the hyperbolic radius as before. If $D$ intersects $\hat{\mathbb{C}} - U$, let $\beta = \beta(D) \in [0, \pi)$ denote the dihedral angle of the intersection. We define $\rho_{\mathrm{hyp}}(D)$ to be the symbol $\infty^\beta$. In particular, $\rho_{\mathrm{hyp}}(D) = \infty^0$ if $D$ is internally tangent to the unit circle $\partial U$. We make the convention that for any angles $\beta_2 \geq \beta_1 \geq 0$ and for any real number $a$, we have: $\infty^{\beta_2} \geq \infty^{\beta_1} > a$. With the same proof, we have:

LEMMA 2.3. *Part* (b) *of Lemma 2.2 still holds if instead requiring the disks of $P$ and $P^*$ be contained in $U$, we assume that all the disks have non-empty intersection with $U$.*

As an application, let us prove Rigidity Theorem 1.2.

*Proof of Theorem 1.2.* Let $P$ and $P^*$ be as given by the theorem. Let $v_0$ be an arbitrary vertex of $G$. We will show that $\rho_{\mathrm{hyp}}(P(v_0)) = \rho_{\mathrm{hyp}}(P^*(v_0))$ and this implies the theorem. By contradiction, let us assume that they are different, say, $\rho_{\mathrm{hyp}}(P(v_0)) < \rho_{\mathrm{hyp}}(P^*(v_0))$. Then for some real $\delta = 1 + \varepsilon > 1$, we still have $\rho_{\mathrm{hyp}}(\delta P(v_0)) < \rho_{\mathrm{hyp}}(P^*(v_0))$. Consider the pattern $\delta P = \{\delta P(v); P(v) \in P\}$. Let $G_1$ be the subgraph of $G$ corresponding to the subpattern $P_1$ of $\delta P$ consisting of those disks which intersect $U$. Let $P_1^*$ be the corresponding subpattern of $P^*$. Then $G_1$ is finite as $P$ is locally finite in $U$. For every boundary vertex $v$ of $G_1$, the disk $P_1(v)$ intersects the boundary of $U$; therefore, $\rho_{\mathrm{hyp}}(P_1(v)) > \rho_{\mathrm{hyp}}(P^*(v))$. Thus, by Lemma 2.3, we have $\rho_{\mathrm{hyp}}(P_1(v_0)) \geq \rho_{\mathrm{hyp}}(P_1^*(v_0))$, a contradiction since $P_1(v_0) = \delta P(v_0)$ and $P_1^*(v_0) = P^*(v_0)$. $\qquad\square$

## 3. Immersed disk patterns

Let $T$ be a triangulation of a planar surface possibly with boundary. Let $G = T^{(1)}$, and let $\Theta : E \to [0, \pi/2]$ be a function on the set of edges. Given a function $\eta : V \to \mathbb{R}_+ = (0, \infty)$ on the set of vertices, we will build a path metric on $|T|$ as follows. For any edge $[v_1, v_2]$, let $D_1$ and $D_2$ be a pair of disks in $\mathbb{C}$ of radius $\eta(v_1)$ and $\eta(v_2)$ respectively, so that their dihedral angle is equal to $\Theta([v_1, v_2])$. Let $\rho([v_1, v_2])$ denote the distance of the centers of $D_1$ and $D_2$. Then there is a unique path metric in $|T|$ in which every edge $[v_1, v_2]$ is isometric to a line segment of length $\rho([v_1, v_2])$, and every 2-simplex is isometric to an euclidean triangle. The curvature $K(v)$ at an interior vertex $v$ is defined to be the sum of the angles at $v$ of the 2-simplexes which contain $v$, less $2\pi$. If $K(v) = 0$ at all interior vertices, then $|T|$ is locally euclidean. If, in addition, $|T|$ is simply connected, then the path metric space has an immersion into the plane. In this case, it is easy to see that there is a collection $P$ of closed disks in the plane, indexed by $V$, such that $\rho(P(v)) = \eta(v)$ for any vertex of $T$ and that the dihedral angle of $P(v_1)$ and $P(v_2)$ is equal to $\Theta([v_1, v_2])$ for any edge $[v_1, v_2]$ in $T^{(1)}$. We will call $P$ an *immersed disk pattern* which realizes the data $(T^{(1)}, \Theta)$. Given $(T^{(1)}, \Theta)$, the isometric class of $P$ is uniquely determined by the radius function $\eta(v) = \rho(P(v))$. The following is clear.

LEMMA 3.1. *Lemmas 2.1, 2.2 and 2.3 also hold for immersed disk patterns.*

Let $\rho_0$, $\rho_1$, $\rho_2$ be positive real numbers, and let $\Theta_{ij} \in [0, \pi/2]$, where $(i, j) \in \{(0, 1), (0, 2), (1, 2)\}$. Then there is a configuration of three disks $D_0$, $D_1$, $D_2$, unique up to euclidean isometries, such that $\rho(D_k) = \rho_k$, $0 \le k \le 2$, and that the dihedral angle of $D_i$ and $D_j$ is equal to $\Theta_{ij}$. The configuration can be shown to exist by first computing (using the rule of cosines) the size of the triangle formed by the centers $A_k$ of $D_k$, $k = 1, 2, 3$ (see e.g. [16]). Let $\Phi_k = \Phi_k(\rho_0, \rho_1, \rho_2; \Theta_{01}, \Theta_{02}, \Theta_{12})$ be the angle of the triangle $A_0 A_1 A_2$ at $A_k$ (see Figure 3.1).

The following lemma is contained in the argument of [26] (or [16, Lemma 3]).

LEMMA 3.2. *Following above, let $L_{ij}$ be the line through $\partial D_i \cap \partial D_j$, and tangent to $D_i$ and $D_j$ in case $\Theta_{ij} = 0$. The lines $L_{01}$, $L_{02}$ and $L_{12}$ meet in a point $O$ contained in the triangle $A_0 A_1 A_2$ (see Figure 3.1). Let $h_{ij}$ denote the distance from $O$ to the edge $A_i A_j$. For any $0 \le i, j \le 2$, $i \ne j$, we have*

$$\frac{\partial \Phi_i(\rho_0, \rho_1, \rho_2; \Theta_{01}, \Theta_{02}, \Theta_{12})}{\partial(\log \rho_j)} = \frac{h_{ij}}{|A_i - A_j|}.$$

*In particular, $\dfrac{\partial \Phi_i}{\partial(\log \rho_j)}$ does not change if one switches $i$ and $j$.*

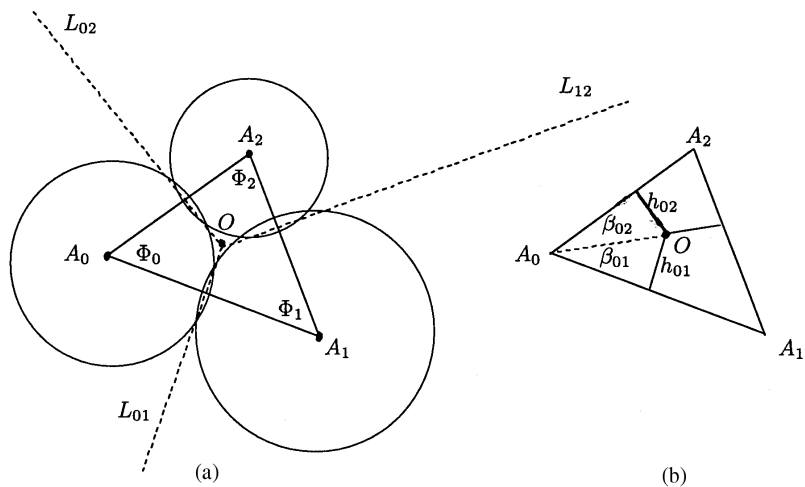

**Figure 3.1.** The configuration of the disks $P(v_i)$ and the corresponding contact graph.

Note that the point $O$ is in the interior of the triangle $A_0 A_1 A_2$ if and only if the angles $\Theta_{01}$, $\Theta_{02}$ and $\Theta_{12}$ are not permutations of $0$, $\pi/2$ and $\pi/2$.

The following elementary lemma will be useful.

LEMMA 3.3.   *Let $j = 1$ or $2$. There is a universal constant $C$ such that:*

$$\frac{\partial \Phi_0(\rho_0, \rho_1, \rho_2; \Theta_{01}, \Theta_{02}, \Theta_{12})}{\partial(\log \rho_j)} \leq C \Phi_0(\rho_0, \rho_1, \rho_2; \Theta_{01}, \Theta_{02}, \Theta_{12}).$$

*Proof.* Without loss of generality, let us assume that $j = 1$. Let $\beta_{0k} = \angle O A_0 A_k$, where $k = 1, 2$. As we remarked earlier, as $\Theta_{01}, \Theta_{02}, \Theta_{12} \leq \pi/2$, point $O$ is contained in the triangle $A_0 A_1 A_2$ (and $O$ is not equal to $A_0, A_1$ or $A_2$). Thus, $\beta_{0j} \in [0, \pi/2)$.

By Lemma 3.2, $\partial \Phi_0/\partial(\log \rho_1) = h_{01}/|A_0 - A_1| < \tan \beta_{01}$, which is less than $C\beta_{01} \leq C\Phi_0$ provided that $\beta_{01}$ is not close to $\pi/2$, say, for $\beta_{01} \leq 1$. Next, consider the case $\beta_{01} \in [1, \pi/2)$. In particular, $\Phi_0 \geq \beta_{01} \geq 1$ is bounded from below. Thus we need show that $\partial \Phi_0/\partial(\log \rho_1) = h_{01}/|A_0 - A_1|$ is bounded by some constant $C$. By rescaling we may assume that $|A_0 - A_1| = 1$. Then $\rho_0, \rho_1 < 1$. But $O$ lies in the convex hull of the union of disks $D_0 \cup D_1$. Thus, $h_{01} < 2$. We conclude that $\partial \Phi_0/\partial(\log \rho_1) = h_{01}/|A_0 - A_1| < 2$.   □

The following lemma is a straight generalization of Andreev's theorem as interpreted by Thurston [26].

Lemma 3.4. *Let $G$ be the one-skeleton of a triangulation $T$ of the closed topological disk, and let $\Theta : E \to [0, \pi/2]$. Assume conditions* (C1) *and* (C2) *hold. For any boundary vertex $v$ of $G$, let $\rho(v) > 0$ be given. Then there is an immersed disk pattern $Q$ in the plane which realizes $(G, \Theta)$, such that $\rho(Q(v)) = \rho(v)$ for any boundary vertex $v$. Moreover, $Q$ is unique up to euclidean isometries.*

*Proof.* We will just sketch the proof here as it is almost identical to the argument of [26] (or [16]). Let $\mathbb{R}_+^N$ denote the space of functions $\mu : V_{\text{int}} \to \mathbb{R}_+ = (0, +\infty)$ where $V_{\text{int}}$ is the set of interior vertices of $G$ and $N$ is the cardinality of $V_{\text{int}}$. Similarly, let $\mathbb{R}^N$ denote the set of functions $K : V_{\text{int}} \to \mathbb{R}$. We construct a transformation $F : \mathbb{R}_+^N \to \mathbb{R}^N$ as follows. For any $\mu$ in $\mathbb{R}_+^k$, define $\eta : V \to \mathbb{R}_+$ by: $\eta(v) = \mu(v)$ for $v \in V_{\text{int}}$ and $\eta(v) = \rho(v)$ for $v \in V - V_{\text{int}}$. Then, as above, $\eta$ defines a path metric in $|T|$. Let $K(v)$, $v \in V_{\text{int}}$, be the curvature at $V$ in the path metric defined by $\eta$. Then our transformation $F$ takes $\mu$ to $K$. By the same argument as in [26], $F$ is one-to-one, and maps $\mathbb{R}_+^N$ onto some open region in $\mathbb{R}^N$ bounded by a finite number of hyperplanes. Using conditions (C1) and (C2), the argument of [16] can be repeated without any essential modification to show that $K = 0 \in \mathbb{R}^N$ is always on the correct side of the hyperplanes, and hence is in the image of $F$. The corresponding preimage $\eta : V_{\text{int}} \to \mathbb{R}$ will then define an immersed disk pattern $Q$ in $\mathbb{C}$ which realizes $(T^{(1)}, \Theta)$ and satisfies the required boundary condition. Clearly, the uniqueness of $Q$ up to euclidean isometries follows by the uniqueness of $\eta$. $\square$

We note that with a similar proof, Lemma 3.4 also holds in the hyperbolic plane.

## 4. Deformation of disk patterns and harmonicity

Unless otherwise specified, we will let $C$ be any positive constant which is independent of the vertices of $G$. The following lemma will be proved in Section 6.

Lemma 4.1. *Let $P$ and $P^*$ be as in Theorem* 1.1. *There is a constant $C \geq 1$ such that for any vertex $v$,*

$$(4.1) \qquad \frac{1}{C} \leq \frac{\rho(P^*(v))}{\rho(P(v))} \leq C.$$

In the following two sections, we will prove the conclusion of Theorem 1.1 assuming Lemma 4.1 by exploring the "harmonicity" of the bounded function $\log \left[ \rho(P^*(v)) / \rho(P(v)) \right]$. Harmonicity properties were also used in [3], [19], [20], [1], and [25] in their study of the rigidity properties of some special classes of disk packings. Our approach may be compared with that of K. Stephenson [25] who studied the convergence of a certain class of disk packings by using random walks and Markov processes in some networks with "leaks" which are based on a graph expanded from $G$. The networks used in this paper are based on $G$, and the edge conductances are quite simple.

Let $G$, $\Theta$, $P$ and $P^*$ be as given in Theorem 1.1. We will first build a one-parameter family of immersed patterns joining $P$ and $P^*$ as follows. For each vertex $v$ in $G$, denote $\tilde{\rho}(v,t) = e^{\lambda(v)t}\rho(P(v))$, where

$$\lambda(v) = \log\big(\rho(P^*(v))/\rho(P(v))\big).$$

Clearly $\tilde{\rho}(v,0) = \rho(P(v))$, $\tilde{\rho}(v,1) = \rho(P^*(v))$, and by (4.1),

$$(4.2) \qquad |\log\frac{\tilde{\rho}(v,t+h)}{\tilde{\rho}(v,t)}| = |h| \cdot |\lambda(v)| \leq |h| \cdot \log C,$$

where $t$, $t + h \in [0,1]$.

Let $G_n$ be an increasing sequence of subgraphs of $G$ whose union is $G$. We will choose $G_n$ to be the one-skeletons of triangulations of the closed topological disk. The pair $(G_n, \Theta|_{G_n})$ is realized by a subpattern of $P$, and thus satisfies conditions (C1) and (C2). Using Lemma 3.4, for each $n$ and each $t \in [0,1]$, there is an immersed disk pattern $P_{n,t}$ in $\mathbb{C}$ which realizes $(G_n, \Theta|_{G_n})$ such that $\rho(P_{n,t}(v))) = \tilde{\rho}(v,t)$ for each boundary vertex $v$ of $G_n$. Using Lemma 2.1 for immersed disk patterns (see Lemma 3.1), it follows from (4.2) that for any vertex $v$ of $G_n$,

$$\left|\log\frac{\rho(P_{n,t+h}(v))}{\rho(P_{n,t}(v))}\right| \leq |h| \cdot \log C.$$

Therefore,

$$(4.3) \qquad \left|\frac{d\log\rho(P_{n,t}(v))}{dt}\right| \leq \log C.$$

In particular, this implies that $\rho(P(v))/C \leq \rho(P_{n,t}(v)) \leq C\rho(P(v))$. Thus, replacing by a subsequence and by means of euclidean transformations if necessary, we may assume that for each vertex $v$ in $G$, the sequence $P_{n,t}(v)$ converges (in the Hausdorff metric) to some limit disk, say $P_t(v)$. Then for each $t$, the collection $P_t = \big\{P_t(v);\ v \in V\big\}$ is an immersed disk pattern.

By uniqueness part of Lemma 3.4, we may assume that the patterns $P_{n,0}$ are subpatterns of $P$, and therefore $P_0 = P$. Similarly, we may assume that $P_1 = P^*$.

Letting $n \to \infty$ in (4.3), we obtain

$$(4.4) \qquad \left|\frac{d\log\rho(P_t(v))}{dt}\right| \leq \log C$$

where the derivative with respect to $t$ is in the generalized sense (or distributional sense). Let $l(v,t) = \log \rho(P_t(v)) - \log \rho(P_0(v))$. We have $l(v,0) = 0$, and (4.4) implies that,

$$(4.5) \qquad \left| \frac{d\,l(v,t)}{dt} \right| \le \log C.$$

Let $t \in [0,1]$ be fixed. Define $h(v) = h(v,t) = d\,l(v,t)/dt$. We will show that $h$ is harmonic in the electrical network based on $G$ in which the conductance $\mu([v_0, v_1])$ of an edge $[v_0, v_1]$ is defined to the partial derivative of the curvature $K(v_0)$ with respect to the logarithm of $\rho(P_t(v_1))$. Let $v_2$ and $v_3$ be the vertices which are neighbors to both vertices $v_0$ and $v_1$ of the edge. That is, $[v_0, v_1, v_2]$ and $[v_0, v_1, v_3]$ are the two 2-simplexes which contain $[v_0, v_1]$. Then we have,

$$(4.6)$$

$$\mu([v_0,v_1]) = \frac{\partial \Phi_0(\rho(P_t(v_0)), \rho_1, \rho(P_t(v_2)); \Theta([v_0,v_1]), \Theta([v_0,v_2]), \Theta([v_1,v_2]))}{\partial \log \rho_1} \bigg|_{\rho_1 = \rho(P_t(v_1))}$$

$$+ \frac{\partial \Phi_0(\rho(P_t(v_0)), \rho_1, \rho(P_t(v_3)); \Theta([v_0,v_1]), \Theta([v_0,v_3]), \Theta([v_1,v_3]))}{\partial \log \rho_1} \bigg|_{\rho_1 = \rho(P_t(v_1))}.$$

Note that $\mu([v_0,v_1]) = \mu([v_1,v_0])$ by Lemma 3.2, so the conductance is well defined. Clearly, $\mu(e) \ge 0$ for any edge $e$, and $\mu(e) = 0$ if and only if $e$ is reducible. Thus, $\mu > 0$ precisely in the reduced graph of $P_t$.

In virtue of Lemma 3.2, the conductance $\mu([v_0,v_1])$ has the following geometrical interpretation. Let $A_0$ and $A_1$ be the centers of $P_t(v_0)$ and $P_t(v_1)$, respectively. Consider the triple $v_0, v_1, v_2$, where there is a unique (real or imaginary) circle which is orthogonal to $\partial P_t(v_0)$, $\partial P_t(v_1)$ and $\partial P_t(v_2)$. Let the center of the circle be denoted by $O_2$; it is called the center of the circle dual to the disks $P_t(v_0)$, $P_t(v_1)$ and $P_t(u)$. It is a good elementary exercise to show that $O_2$ is at the intersection of the lines $L_{01}, L_{02}$ and $L_{12}$ in Figure 3.1. Similarly, let $O_3$ be the center of the center dual to the disks $P_t(v_0)$, $P_t(v_1)$ and $P_t(v_3)$. Clearly, the line segment $O_2O_3$ is orthogonal to $A_0A_1$. Then the conductance $\mu([v_0,v_1])$ of $[v_0, v_1]$ is equal to the ratio of the euclidean lengths of $O_2O_3$ and $A_0A_1$.

LEMMA 4.2. *The function $h(v) = h(v,t)$ is harmonic in the network based on $G$ in which the conductance $\mu([v_0,v_1])$ of an edge $[v_0, v_1]$ is defined by (4.6).*

*Proof.* Let $v_0$ be an arbitrary vertex and let $v_1, v_2, \ldots, v_l, v_{l+1} = v_1$ be the chain of neighboring vertices which surround $v_0$. Then the Laplacian $\Delta h$ is defined to be

$$(4.7) \qquad \Delta h(v_0) = \sum_{k=1}^{l} \mu([v_0, v_k])(h(v_k) - h(v_0)).$$

Consider the pattern $P_t$. Let $A_k(t)$ denote the center of the disk $P_t(v_k)$, $0 \le k \le l+1$. For $1 \le k \le l$, denote,

$$\begin{aligned}
\beta_k(t) &= \angle A_k(t) A_0(t) A_{k+1}(t) \\
&= \Phi_0\Big(\rho(P_t(v_0)), \rho(P_t(v_k)), \rho(P_t(v_{k+1})); \Theta([v_0,v_k]), \Theta([v_0,v_{k+1}]), \Theta([v_k,v_{k+1}])\Big).
\end{aligned}$$

Then the curvature at $v_0$ is $K(v_0, P_t) = \sum_{1 \le k \le l} \beta_k(t) - 2\pi = 0$ (see Figure 4.1).

Since $(G, \Theta)$ is fixed, $K(v_0, P_t)$ can be considered as a function of the radii $\rho(P_t(v_j))$, $0 \le j \le l$. By the definition of $\mu$, we have,

$$\partial K(v_0, P_t)/\partial \log(\rho(P_t(v_k))) = \mu([v_0, v_k]), \quad 1 \le k \le l.$$

Clearly, $K(v_0, P_t)$ is homogeneous of degree 0 in the (vector) variable $(\rho(P_t(v_0)), \rho(P_t(v_1)), \dots, \rho(P_t(v_l)))$; thus, $\sum_{0 \le j \le l} \partial K(v_0, P_t)/\partial \log(\rho(P_t(v_j))) = 0$. Then, $\partial K(v_0, P_t)/\partial \log(\rho(P_t(v_0))) = -\sum_{1 \le k \le l} \mu([v_0, v_k])$. Differentiating the equality $K(v_0; P_t) = 0$ with respect to $t$, we obtain $\Delta h(v_0) = 0$ by the chain rule. $\qquad\square$

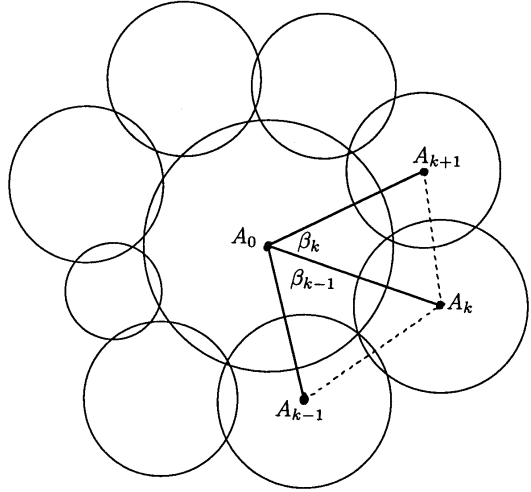

**Figure 4.1.** The curvature at $v_0$ is $\sum_{1 \le k \le l} \beta_k - 2\pi$.

The following lemma follows easily by Lemma 3.3.

LEMMA 4.3.   *There is a universal constant $C$ such that for any vertex $v$ of $G$, the sum of conductances $\mu$ of edges from $v$ is bounded by $C$.*

In the next section, we will show that the network based on $G$ in which the edge conductances are equal to $\mu$ is recurrent, and as a consequence, any bounded harmonic function in the network is a constant.

## 5. Recurrence of electrical networks

We will make extensive use of the vertex extremal length in graphs. The notion resembles the extremal length of curve families in Riemann surfaces, and was introduced and studied by J. Cannon and others (see [6], [7], [23], and [13]). (A further extension to the idea of "transboundary extremal length" appeared in [24].) Below we will recall the basics. For details, we refer to [13].

Let $G$ be a graph and let $V$ be the set of vertices. A *vertex metric* in the graph is just a function $\eta : V \to [0, +\infty)$. The *area* of the metric $\eta$ is defined to be:

$$(5.1) \qquad \text{area}\,(\eta) = \sum_{v \in V} \eta(v)^2.$$

Any subset of vertices will be called a *vertex curve*. The *length* of a vertex curve $\gamma$ in the metric $\eta$ is defined to be:

$$(5.2) \qquad \int_\gamma d\eta = \sum_{v \in \gamma} \eta(v).$$

Let $\Gamma$ be a collection of vertex curves in $G$. A vertex metric $\eta$ is called $\Gamma$-*admissible*, if $\int_\gamma d\eta \geq 1$ for each $\gamma \in \Gamma$. The *vertex modulus* of $\Gamma$ is defined by,

$$(5.3) \qquad \text{MOD}\,(\Gamma) = \inf \big\{ \text{area}\,(\eta); \ \eta : V \to [0, \infty) \text{ is } \Gamma\text{-admissible} \big\}.$$

The *vertex extremal length* of $\Gamma$ is defined to be:

$$(5.4) \qquad \text{VEL}\,(\Gamma) = \frac{1}{\text{MOD}\,(\Gamma)}.$$

We make the convention that the vertex extremal length is $+\infty$ if $\Gamma$ is void.

A *path* in $G$ is defined to be a finite or infinite sequence of vertices $(v_0, v_1, \ldots, v_k, \ldots)$, such that $[v_{k-1}, v_k]$ is an edge of $G$. The set of vertices in a path is a vertex curve, and for the purpose of defining vertex extremal length, we will identify the path with the curve. (Later we will also consider a path as a sequence of edges for the purpose of defining the electrical resistance.) Let $V_1$, $V_2$ be nonvoid subsets of vertices. Let $\Gamma_G(V_1, V_2)$ be the set of paths in $G$ joining $V_1$ and $V_2$. We allow $V_2 = \infty$, in which case a path joining $V_1$ and $\infty$ is by definition a path starting from a vertex in $V_1$ which passes through an infinite subset of vertices. The *vertex extremal length* between $V_1$ and $V_2$ in $G$ is then defined to be:

$$(5.5) \qquad \text{VEL}\,(V_1, V_2) = \text{VEL}\,(\Gamma_G(V_1, V_2)).$$

For any three nonvoid subsets of vertices $V_1$, $V_2$ and $V_3$ in the graph $G$, $V_2$ is said to *separate* $V_1$ and $V_3$, if any path from $V_1$ to $V_3$ passes through a vertex in $V_2$. Note that we *do not* require that $V_i \cap V_j = \emptyset$ for $1 \le i < j \le 3$. In graph theory, $V_2$ is also called a *cutset* between $V_1$ and $V_3$. The lemma below follows by a classical argument (see e.g. [2] or [15]).

LEMMA 5.1. *Let* $V_1$, $V_2$, ..., $V_{2m}$ *be mutually disjoint, nonvoid subsets of vertices such that for* $i_1 < i_2 < i_3$, $V_{i_2}$ *separates* $V_{i_1}$ *from* $V_{i_3}$. *We allow* $V_{2m} = \infty$. *Then we have*

$$(5.6) \qquad \mathrm{VEL}\,(V_1, V_{2m}) \ge \sum_{i=1}^{m} \mathrm{VEL}\,(V_{2k-1}, V_{2k}).$$

Let $V_0$ be a finite, nonvoid subset of vertices in a *connected* graph $G$. Then $G$ is called VEL-*parabolic* if $\mathrm{VEL}\,(V_0, \infty) = \infty$, and VEL-*hyperbolic* otherwise. The definition is independent of the choice of $V_0$.

If $G$ is the contact graph of a disk packing which is locally finite in $\mathbb{C}$, and if $G$ is connected, then $G$ is VEL-parabolic (see [13, Theorem 3.1]). Without much more effort, we will prove:

LEMMA 5.2. *Let* $G$ *be a connected graph, and let* $\Theta : E \to [0, \pi/2]$. *Suppose that* $P$ *is a disk pattern in* $\mathbb{C}$ *which realizes* $(G, \Theta)$. *If* $P$ *is locally finite in the plane, then* $G$ *is* VEL-*parabolic.*

*Proof.* For a circle $c$ in the plane, let $V_c = V_c(P) = \{v;\ P(v) \cap c \ne \emptyset\}$. Then, inductively, we may find a sequence of circles $c(r_i) = \{z \in \mathbb{C};\ |z| = r_i\}$, such that: $r_{i+1} \ge 2r_i$ and $V_{c(r_{i+1})} \cap V_{c(r_i)} = \emptyset$ (and $V_{c(r_1)} \ne \emptyset$). For $1 \le i_1 < i_2 < i_3$, it is easy to see that $V_{c(r_{i_2})}$ separates $V_{c(r_{i_1})}$ from $V_{c(r_{i_3})}$. By Lemma 5.1, we have:

$$\mathrm{VEL}\,(V_{c(r_1)}, \infty) \ge \sum_{k=1}^{\infty} \mathrm{VEL}\,(V_{c(r_{2k-1})}, V_{c(r_{2k})}).$$

Then Lemma 5.2 follows by the lemma below. $\qquad \square$

LEMMA 5.3. *Let* $G$ *be a graph, and let* $\Theta : E \to [0, \pi/2]$. *Suppose that* $P$ *is a disk pattern in* $\mathbb{C}$ *which realizes* $(G, \Theta)$. *Then for any* $r_2 > r_1 > 0$,

$$(5.7) \qquad \mathrm{VEL}\,(V_{c(r_1)}, V_{c(r_2)}) \ge \frac{(r_2 - r_1)^2}{(32 + (8\pi)^2)(r_2)^2}.$$

*In particular, if* $r_2 \ge 2r_1$, *then,*

$$(5.8) \qquad \mathrm{VEL}\,(V_{c(r_1)}, V_{c(r_2)}) \ge \frac{1}{128 + (16\pi)^2}.$$

*Remark.* We first note that the constants in (5.7), (5.8), and in similar formulas below are not the best possible. In fact, any constants, large or small, are good for our purpose.

*Proof.* For any $v \in V$, let $d(v) = \operatorname{diam}(P(v) \cap D(r_2))$, where $\operatorname{diam}(\cdot)$ denotes the euclidean diameter. Let $\gamma$ be a path in $G$ joining $V_{c(r_1)}$ and $V_{c(r_2)}$. Then $P(\gamma) = \bigcup_{w \in \gamma} P(w)$ is a continuum joining $c(r_1)$ and $c(r_2)$ (see Figure 5.1).

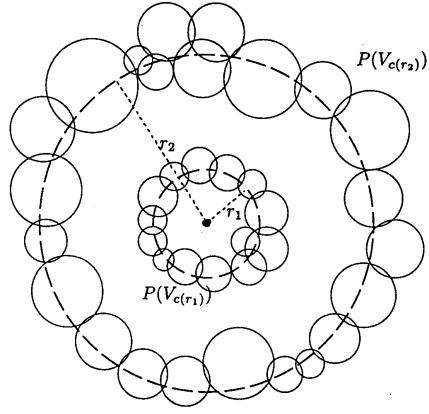

**Figure 5.1.** For any path $\gamma$ in $G$ joining $V_{c(r_1)}$ and $V_{c(r_2)}$, $P(\gamma)$ is a continuum joining $c(r_1)$ and $c(r_2)$.

It follows that $\sum_{v \in \gamma} d(v) \geq r_2 - r_1$. Let $\eta(v) = d(v)/(r_2 - r_1)$. Thus $\eta$ is $\Gamma_G(V_{c(r_1)}, V_{c(r_2)})$-admissible, and then by (5.3),

$$
\begin{aligned}
\operatorname{MOD}\left(\Gamma_G(V_{c(r_1)}, V_{c(r_2)})\right) &\leq \operatorname{area}(\eta) \\
&\leq \sum_{\operatorname{area}(P(v) \cap D(r_2)) \geq (1/2)\operatorname{area}(P(v))} (4/\pi)\operatorname{area}(P(v))/(r_2 - r_1)^2 \\
&+ \sum_{\operatorname{area}(P(v) \cap D(r_2)) < (1/2)\operatorname{area}(P(v))} d(v)^2/(r_2 - r_1)^2.
\end{aligned}
$$

Since any point in $\mathbb{C}$ lies in at most four disks in the pattern $P$, the first sum on the right-hand side is bounded by $(4/\pi) \cdot 4 \cdot 2 \cdot \operatorname{area}(D(r_2))/(r_2 - r_1)^2 = 32(r_2)^2/(r_2 - r_1)^2$. The second sum is bounded by $(8\pi r_2)^2/(r_2 - r_1)^2$, because,

$$
\sum_{\operatorname{area}(P(v) \cap D(r_2)) < (1/2)\operatorname{area}(P(v))} d(v) \leq 4 \cdot \operatorname{length}(\partial D(r_2)) = 8\pi r_2.
$$

As a consequence,

$$
(5.9) \qquad \operatorname{MOD}\left(\Gamma_G(V_{c(r_1)}, V_{c(r_2)})\right) \leq \operatorname{area}(\eta) \leq [32 + (8\pi)^2](r_2)^2/(r_2 - r_1)^2.
$$

This implies (5.7).                                                                 $\square$

Let $\mu: E \to [0, \infty)$ be a function on the set of edges in a graph $G$. Let us denote by $G_\mu$ the electrical network based on $G$ in which an edge $e$ is placed with a resistor of conductance $\mu(e)$ (see e.g. [8]). The *reduced graph* of the network $G_\mu$ is the graph formed by the same vertex set as $G$ and those edges $e$ of $G$ for which $\mu(e) \neq 0$. A network is called *connected* if its reduced graph is connected.

An *edge metric* in $G$ is by definition a function $m : E \to [0, \infty)$. Its *area* (or $\mu$-area) is defined by,

$$(5.10) \qquad \operatorname{area}_\mu(m) = \sum_{e \in E} \mu(e) m(e)^2.$$

The *length* of a path $\gamma = (v_0, v_1, \ldots, v_k, \ldots)$ in the metric $m$ is

$$(5.11) \qquad \int_\gamma dm = \sum_{k \geq 1} m([v_{k-1}, v_k]).$$

Let $\Gamma$ be a collection of paths. We call the edge metric $\Gamma$-*admissible* if $\int_\gamma dm \geq 1$ for any $\gamma \in \Gamma$. The *conductance* (or $\mu$-*conductance*) of $\Gamma$ is defined by,

$$(5.12) \qquad \operatorname{COND}(\Gamma) = \inf \left\{ \operatorname{area}_\mu(m);\ m : E \to [0, \infty) \text{ is } \Gamma\text{-admissible} \right\}.$$

Its inverse is called the *resistance* (or $\mu$-*resistance*) of $\Gamma$,

$$(5.13) \qquad \operatorname{RES}(\Gamma) = \frac{1}{\operatorname{COND}(\Gamma)}.$$

In case $\Gamma = \Gamma_G(V_1, V_2)$ where $V_1$ and $V_2$ are a pair of mutually disjoint, nonvoid subsets of vertices, we call $\operatorname{COND}(\Gamma)$ and $\operatorname{RES}(\Gamma)$ the *conductance* and *resistance between* $V_1$ and $V_2$. We denote $\operatorname{COND}(V_1, V_2) = \operatorname{COND}(\Gamma_G(V_1, V_2))$, and $\operatorname{RES}(V_1, V_2) = \operatorname{RES}(\Gamma_G(V_1, V_2))$. Again, we allow $V_2 = \infty$. A connected network $G_\mu$ is called *recurrent* (or *parabolic*) if the resistance from a finite nonvoid vertex subset to the infinity is infinite; and *transient* (or *hyperbolic*) otherwise.

LEMMA 5.4. *Let $C > 0$ be a constant. Suppose that for each vertex $v$, we have $\sum_{v \in e} \mu(e) \leq C$. Then for any mutually disjoint, nonvoid subsets $V_1$ and $V_2$ of $V$, we have:*

$$(5.14) \qquad \operatorname{VEL}(V_1, V_2) \leq 2C \cdot \operatorname{RES}(V_1, V_2).$$

*In particular, if $G$ is* VEL-*parabolic and if $G_\mu$ is connected, then $G_\mu$ is recurrent.*

*Proof.* We need show that,

$$\operatorname{COND}(V_1, V_2) \leq 2C \cdot \operatorname{MOD}(V_1, V_2).$$

Let $\eta : V \to [0, \infty)$ be any $\Gamma_G(V_1, V_2)$-admissible vertex metric. Define $m = m_\eta : E \to [0, \infty)$ by, $m([v_1, v_2]) = \eta(v_1) + \eta(v_2)$. Then clearly, $m$ is a $\Gamma_G(V_1, V_2)$-admissible edge metric. The lemma follows since,

$$\text{area}_\mu(m) = \sum_{e=[v_1,v_2]\in E} \mu(e)(\eta(v_1)+\eta(v_2))^2 \le 2\sum_{e=[v_1,v_2]\in E} \mu(e)\big(\eta(v_1)^2+\eta(v_2)^2\big)$$
$$= 2\sum_{v\in V}\big(\sum_{v\in e}\mu(e)\big)\eta(v)^2 \le 2C\sum_{v\in V}\eta(v)^2 = 2C\cdot\text{area}(\eta). \qquad \square$$

LEMMA 5.5. *Let $G_\mu$ be a recurrent network, let $V$ be the set of vertices, and let $h : V \to \mathbb{R}$ be a bounded harmonic function. Then $h$ is constant.*

*Proof.* The lemma is well known to analysts. We will sketch an elementary proof for those who are not familiar with the method. For a survey of relevant results and their origin, we refer to [27].

First, we may assume that $\mu > 0$ on all edges, otherwise replace $G$ by the reduced graph. A vertex $v$ of a subset $W$ of vertices is called an *interior vertex* of the subset *relative to $G$*, if its neighbors in $G$ are all contained in the subset; otherwise it is called a *boundary vertex* of $W$ *relative to $G$*. Let $W$ be finite and let $g : W \to \mathbb{R}$ be a function. If for each interior vertex $v$ of $W$ relative to $G$,

$$\Delta g(v) = \sum_{e=[v,w]\in E} \mu(e)(g(w)-g(v)) = 0,$$

then the maximum (or minimum) of $g$ is obviously attained at some boundary vertex of $W$ relative to $G$. This property is referred to as the maximum principle.

Let $v_0$ be fixed. We will show that $h(v) = h(v_0)$ for all $v$. We may reduce to the case when $h(v_0) = 0$ and $|h(v)| \le 1$, so we make this assumption. For each integer $n \ge 0$, let $U_n$ be the set of the vertices which can be joined to $v_0$ by a path of combinatorial length $\le n$. In particular, $U_0 = \{v_0\}$. Let $\partial U_n = U_n - U_{n-1}$, $n \ge 1$. Then by an elementary argument,

$$\lim_{n\to\infty} \text{COND}\,(U_0, \partial U_n) = \text{COND}\,(U_0, \infty) = 0. \tag{5.15}$$

On the other hand, it is easy to check that,
$$\tag{5.16}$$
$$\text{COND}\,(U_0, \partial U_n) = \inf\big\{\text{area}_\mu\,(|\bigtriangledown g|);\ g : V \to \mathbb{R},\ g|_{U_0} = 0,\ g|_{\partial U_n} = 1\big\},$$

where $|\bigtriangledown g|([v,w]) = |g(v) - g(w)| : E \to \mathbb{R}$ is the absolute gradient of $g$. Let $g_n : V \to [0,1]$ be the minimizer for (5.16), which exists as $U_n$ is finite. Then for $v \in V - (\partial U_n \cup U_0)$, we have $\Delta g_n(v) = 0$, and as $h$ is harmonic, $\Delta(g_n + h)(v) = 0$. But $g_n + h = 1 + h \ge 1 - 1 = 0$ in $\partial U_n$ and $g_n + h = 0$ in $U_0$. Then by the maximum principle, the minimum of $g_n + h$ in the finite graph $U_n$ is equal to 0. That is,

$$g_n(v) + h(v) \ge 0, \quad \text{for all}\ v \in U_n. \tag{5.17}$$

On the other hand,

$$\sum_{[v_1,v_2]\in E} \mu([v_1,v_2])|g_n(v_1) - g_n(v_2)|^2 = \text{COND}\,(U_0, \partial U_n) \to 0,$$

by (5.15). As $\mu([v_1, v_2]) > 0$ for any $[v_1, v_2] \in E$, we deduce that

$$\lim_{n\to\infty} |g_n(v_1) - g_n(v_2)| = 0.$$

Since the graph is connected and $g_n(v_0) = 0$, we obtain, $\lim_{n\to\infty} g_n(v) = 0$ for any $v$. Taking limit in (5.17), it follows that $h(v) \geq 0$. Similarly, $-h(v) \geq 0$; hence $h(v) = 0$. $\qquad\square$

We now assemble the proof of Theorem 1.1.

*Proof of Theorem* 1.1 (*assuming Lemma* 4.1). Let $t \in [0, 1]$ be fixed. Consider the network $G_\mu$ where $\mu$ is defined by (4.6). Its reduced graph is the same as the reduced graph of the pattern $P$ (see the remark before Lemma 4.2), and hence is connected. Lemma 4.3 says there is a universal upper bound on the sum of the conductances of the edges from a vertex; and Lemma 5.2 tells us that $G$ is VEL-parabolic. Then using Lemma 5.4, we deduce that the network $G_\mu$ is recurrent.

By Lemma 4.2, $h(v, t) = dl(v, t)/dt$ is harmonic in $G_\mu$. By (4.5), it is also uniformly bounded. Then by Lemma 5.5, $h(v, t)$ is independent of $v$. It follows that $l(v, t) = \int_0^t h(v, s)ds$ is also independent of $v$. In particular, $\rho(P^*(v))/\rho(P(v)) = \rho(P_1(v))/\rho(P_0(v)) = e^{l(v,1)}$ is constant. This implies that $P^*$ and $P$ are images of each other by euclidean similarities. $\qquad\square$

## 6. Uniform bound on the ratio of radii

In this section, we prove Lemma 4.1, restated as follows:

LEMMA 6.1. *Let $P$ and $P^*$ be as in Theorem* 1.1. *There is a constant $C \geq 1$ such that for any vertex $v$,*

$$\frac{1}{C} \leq \frac{\rho(P^*(v))}{\rho(P(v))} \leq C.$$

Let $G$ be a graph, and let $V_1$ and $V_2$ be nonvoid subsets of vertices in $G$. Again, we allow $V_2 = \infty$. Let $\Gamma^*_G(V_1, V_2)$ be the collection of vertex curves in $G$ which separate $V_1$ from $V_2$. Then by an elementary argument similar to the duality argument for the extremal length of curve families in the plane (see e.g. [23] or [13, §5]; compare [7]), we have,

$$(6.1) \qquad \text{MOD}\,(\Gamma^*_G(V_1, V_2)) = \frac{1}{\text{MOD}\,(\Gamma_G(V_1, V_2))} = \text{VEL}\,(V_1, V_2).$$

Let $Q$ be a disk pattern in $\mathbb{C}$ which realizes $(G, \Theta)$ where $0 \leq \Theta \leq \pi/2$. If $G$ is a disk triangulation graph and $Q$ is locally finite in $\mathbb{C}$, then for any circle $c$, $V_c = V_c(Q)$ is a connected set of vertices, and thus $Q(V_c) = \bigcup_{v \in V_c} Q(v)$ is pathwise connected. It is also easy to see that for any path $\gamma = (v_0, v_1, \ldots, v_l)$, the set $Q(\gamma) = \bigcup_{k=0}^{l} Q(v_k)$ is pathwise connected.

It is well known that for an annulus of sufficiently big modulus in $\mathbb{C}$ which separates $0$ from $\infty$, there is some $r > 0$ such that the annulus contains $\{z; \ r \leq z \leq 2r\}$ (see e.g. [2] or [15]). Similarly, we have the following lemma.

LEMMA 6.2. *Let $G$ be a disk triangulation graph and $\Theta : E \to [0, \pi/2]$. Let $Q$ be a disk pattern in $\mathbb{C}$ which realizes the data $(G, \Theta)$. Let $V_0 = \{v_0\}$, $V_1$ and $V_2$ be mutually disjoint, finite, connected subsets of vertices, and let $V_4 = \infty$. Assume that for any $0 \leq i_1 < i_2 < i_3 \leq 4$, the set $V_{i_2}$ separates $V_{i_1}$ from $V_{i_3}$. If*

$$\text{(6.2)} \qquad \text{VEL}\,(V_1, V_2) > 288 + (24\pi)^2,$$

*then there is some $r > 0$ such that for any circle $c$ concentric with $Q(v_0)$ and of radius $\rho(c) \in [r, 2r]$, the vertex subset $V_c(Q)$ separates $V_1$ from $V_2$.*

*Proof.* We may assume that $Q(v_0)$ is centered at $0$. Let $r > 0$ be the minimum number such that $Q(v) \cap D(r) \neq \emptyset$ for any $v \in V_1$. In fact, $r = \max\{\text{dist}\,(0, P(v)); \ v \in V_1\}$; where $\text{dist}\,(\cdot, \cdot)$ denotes the euclidean distance (see Figure 6.1). Without loss of generality, we may assume that $r = 1$.

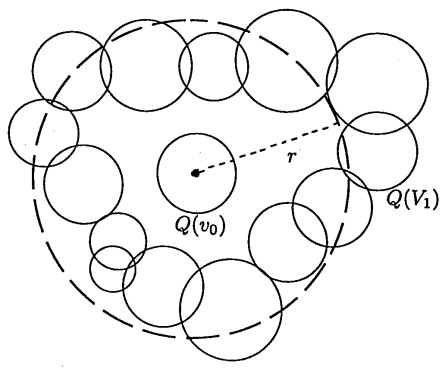

**Figure 6.1.** The number $r$ is the maximum of $\text{dist}\,(0, P(v))$, $v \in V_1$.

Let $\gamma^* \in \Gamma_G^*(V_1, V_2)$. Then $\gamma^*$ also separates $V_1 \cup \{v_0\}$ from $\infty$. Thus by the minimality of $r$ and the assumptions on $V_1$, $V_2$ and $G$, we have,

$$\text{(6.3)} \qquad \text{diam}\,(Q(\gamma^*)) \geq r = 1.$$

By contradiction, let us assume that Lemma 6.2 fails. Then there is a circle $c(\rho_1)$ centered at $z_0 = 0$ and of radius $\rho_1 \in [1, 2]$ such that $V_{c(\rho_1)}$ does not separate $V_1$ from $V_2$. So there is a path $\gamma_0$ in the graph joining $V_1$ and $V_2$ such that $\gamma_0 \cap V_{c(\rho_1)} = \emptyset$. That is, $Q(\gamma_0) \cap c(\rho_1) = \emptyset$. In other words, the continuum $Q(\gamma_0)$ is either contained in the interior of the disk $D(\rho_1)$ or in the exterior of $D(\rho_1)$. But since $D(1)$ has nonempty intersection with every $Q(v)$, $v \in V_1$, and since $\gamma_0$ contains a vertex of $V_1$, we have $Q(\gamma_0) \cap D(1) \neq \emptyset$.

Therefore $Q(\gamma_0)$ is contained in the interior of $D(\rho_1)$ (recall that $\rho_1 \geq 1$). As $\gamma_0 \in \Gamma_G(V_1, V_2)$ and $\gamma^* \in \Gamma_G^*(V_1, V_2)$, we deduce that,

$$(6.4) \qquad Q(\gamma^*) \cap D(\rho_1) \supseteq Q(\gamma^*) \cap Q(\gamma_0) \supseteq Q(\gamma^* \cap \gamma_0) \neq \emptyset.$$

For any vertex $v$, let $\eta(v) = \operatorname{diam}(Q(v) \cap D(3))$. Then $\eta : V \to [0, \infty)$ is a vertex metric in the graph. If $\gamma^* \in \Gamma_G^*(V_1, V_2)$ is *connected*, then by (6.4) it follows that $Q(\gamma^*)$ is either contained in $D(3)$, or is a continuum joining $c(\rho_1)$ and $c(3)$. In either case, by (6.3) and the inequality $3 \geq \rho_1 + 1$, we have, $\int_{\gamma^*} d\eta = \sum_{v \in \gamma^*} \operatorname{diam}(Q(v) \cap D(3)) \geq 1$. Again by the assumptions on $V_1$, $V_2$ and $G$, it is easy to see that each vertex curve in $\Gamma_G^*(V_1, V_2)$ contains a *connected* vertex subcurve $\gamma^*$ in $\Gamma_G^*(V_1, V_2)$. It follows that $\eta$ is $\Gamma_G^*(V_1, V_2)$-admissible. Thus,

$$\operatorname{MOD}(\Gamma_G^*(V_1, V_2)) \leq \operatorname{area}(\eta).$$

By an argument as in Lemma 5.3 (see (5.9)), we deduce that,

$$\operatorname{MOD}(\Gamma_G^*(V_1, V_2)) \leq \operatorname{area}(\eta) \leq [32 + (8\pi)^2] \cdot 3^2.$$

By (6.1), this contradicts (6.2). $\qquad\blacksquare$

**COROLLARY 6.2.** *Let $G$ be a disk triangulation graph and let $\Theta : E \to [0, \pi/2]$ be a function on the set of edges. Let $Q$ be a disk pattern in $\mathbb{C}$ which realizes the data $(G, \Theta)$. Then $Q$ is locally finite in $\mathbb{C}$ if and only if $G$ is VEL-parabolic.*

*Proof.* The "only if" part follows by Lemma 5.2. To prove the converse, let $G$ be VEL-parabolic. Let $v_0$ be a vertex of $G$, and let $V_0 = \{v_0\}$, $V_i$, $i = 1, 2, \ldots$, be a sequence of mutually disjoint, finite, connected subsets of vertices such that for any $0 \leq i_1 < i_2 < i_3$, the set $V_{i_2}$ separates $V_{i_1}$ from $V_{i_3}$, and that,

$$\operatorname{VEL}(V_i, V_{i+1}) > 288 + (24\pi)^2.$$

We may assume that $Q(v_0)$ is centered at 0. Then by Lemma 6.1, for any integer $k \geq 1$, there is some $r_k > 0$, such that for any $\rho \in [r_k, 2r_k]$, the vertex set $V_{c(\rho)}(Q)$ separates $V_{2k-1}$ from $V_{2k}$.

On one hand, this implies that $V_{c(2r_k)}(Q) \cap V_{c(r_{k+1})}(Q) \subseteq V_{2k} \cap V_{2k+1} = \emptyset$, and therefore, $r_{k+1} \geq 2r_k$; hence, $\lim_{k \to \infty} r_k = +\infty$. On the other hand, we deduce that the finite, connected subset $V_{2k}$ separates $V_{c(r_k)}(Q)$ from $\infty$. Thus there are only a finite number of disks in the pattern which intersect the disk $D(r_k)$. It follows that $Q$ is locally finite in the plane. $\qquad\blacksquare$

For any disk $D$ in $\mathbb{C}$, denote,

$$(6.5) \qquad \tau(D) = \frac{\rho(D)}{\operatorname{dist}(0, D)} \in (0, +\infty].$$

For an infinite disk pattern $Q$ in $\mathbb{C}$, we may arrange the disks of $Q$ into a sequence, $Q(v_k)$, and define,

$$(6.6) \qquad \tau(Q) = \limsup_{v \to \infty} \tau(Q(v)) = \limsup_{k \to \infty} \tau(Q(v_k)) \in [0, +\infty].$$

Clearly, $\tau(Q)$ does not depend on the order of the sequence. Moreover, if $Q$ is locally finite in $\mathbb{C}$, and if $f : \mathbb{C} \to \mathbb{C}$ is a euclidean similarity, then $\tau(f(Q)) = \tau(Q)$.

*Proof of Lemma* 6.1. Let $G$, $\Theta$, $P$ and $P^*$ be as in Theorem 1.1. We may assume that for some fixed vertex $v_0$, we have $P^*(v_0) = P(v_0)$ and its center is at 0. Let $v_1$ be any other vertex. For any circle $c$ of $\hat{\mathbb{C}}$, we will denote $V_c = V_c(P)$, and $V_c^* = V_c(P^*)$.

*Case* 1: $\tau(P) < \infty$.   There is some $\delta > 0$ and a finite subset $V_0$ of vertices such that,

$$(6.7) \qquad \frac{\rho(P(v))}{\operatorname{dist}(0, P(v))} = \tau(P(v)) \leq \frac{\delta}{3}, \quad \text{for all} \ \ v \in V - V_0.$$

We will assume that $V_0$ contains the vertices $v_0$ and $v_1$, and that

$$(6.8) \qquad \delta > 3.$$

Let $r_0 > 0$ be such that $P(V_0) \subseteq D(r_0)$. Then for any $r > r_0$, it follows by (6.7) and (6.8) that no disk $P(v)$ can intersect both $c(r)$ and $c(\delta r)$. In other words, $V_{c(r)} \cap V_{c(\delta r)} = \emptyset$.

Let $r_i = \delta^i r_0$, and let $V_i = V_{c(r_i)}$ for $i \geq 1$. Then for any $0 \leq i_1 < i_2 < i_3$, $V_{i_2}$ separates $V_{i_1}$ from $V_{i_3}$. Therefore,

$$\operatorname{VEL}(V_2, V_{2k-1}) \geq \sum_{i=1}^{k-1} \operatorname{VEL}(V_{2i}, V_{2i+1}).$$

By Lemma 5.3, $\operatorname{VEL}(V_{2i}, V_{2i+1}) \geq 1/\big(128 + (16\pi)^2\big)$. It follows that,

$$\operatorname{VEL}(V_2, V_{2k-1}) \geq \frac{k-1}{128 + (16\pi)^2}.$$

Let $k$ be the following integer:

$$(6.9) \qquad k = 2 + \text{the integer part of } (128 + (16\pi)^2)(288 + (24\pi)^2).$$

Then,

$$(6.10) \qquad \operatorname{VEL}(V_2, V_{2k-1}) > 288 + (24\pi)^2.$$

Applying Lemma 6.1 to the pattern $Q = P^*$, there is some $r^*$, such that for any $\rho \in [r^*, 2r^*]$, the vertex set $V_{c(\rho)}^*$ separates $V_2$ from $V_{2k-1}$. It follows that $V_{c(\rho)}^*$ also separates $V_1$ from $V_{2k}$. Moreover, $V_{c(\rho)}^* \cap V_1 \subseteq V_2 \cap V_1 = \emptyset$ and $V_{c(\rho)}^* \cap V_{2k} \subseteq V_{2k-1} \cap V_{2k} = \emptyset$. Therefore $P^*(V_1)$ is contained in the interior disk of $c(r^*)$, and $P^*(V_{2k})$ lies in the exterior disk of $c(2r^*)$. (See Figure 6.2.)

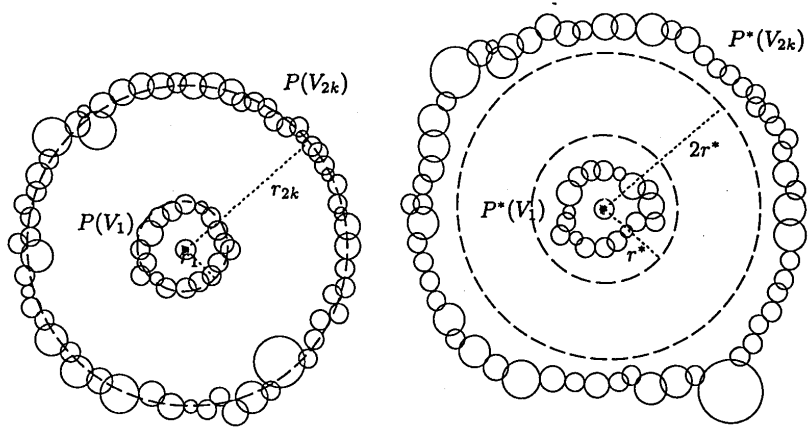

**Figure 6.2.** The sets $P(V_1)$, $P(V_{2k})$, $P^*(V_1)$ and $P^*(V_{2k})$ (not drawn to scale).

Consider the subgraph $G_1$ of $G$ consisting of the vertices $v$ for which $P(v) \cap D(r_1) \neq \emptyset$ and the corresponding edges. Then the disks $P^*(v)$, $v \in V(G_1)$, are contained in the disk $D(r^*) \subset D(2r^*)$. Using Lemma 2.3, we deduce that for each $v \in V(G_1)$,

$$\rho_{\text{hyp}}\left((1/(2r^*))P^*(v)\right) \leq \rho_{\text{hyp}}\left((1/r_1)P(v)\right).$$

In particular, for $v = v_j$, $j = 0, 1$,

$$(6.11) \qquad \rho_{\text{hyp}}\left((1/(2r^*))P^*(v_j)\right) \leq \rho_{\text{hyp}}\left((1/r_1)P(v_j)\right).$$

On the other hand, consider the subgraph $G_2$ consisting of the vertices $v$ for which $P^*(v) \cap D(2r^*) \neq \emptyset$ and the corresponding edges. Then the disks $P(v)$, $v \in G_2$, are contained in the interior of $D(r_{2k})$. If follows by Lemma 2.3 again that,

$$\rho_{\text{hyp}}\left((1/(2r^*))P^*(v)\right) \geq \rho_{\text{hyp}}\left((1/r_{2k})P(v)\right), \quad \text{for all} \;\; v \in V(G_2).$$

Thus, for $v = v_j$, $j = 0, 1$,

$$(6.12) \qquad \rho_{\text{hyp}}\left((1/(2r^*))P^*(v_j)\right) \geq \rho_{\text{hyp}}\left((1/r_{2k})P(v_j)\right).$$

Since the disks $(1/(2r^*))P^*(v_j)$, $(1/r_1)P(v_j)$, $(1/r_{2k})P(v_j)$, $j = 0, 1$, are all contained in $(1/2)U$, it follows that their euclidean radii are comparable with their hyperbolic radii. So by (6.11) and (6.12), there is a universal constant $C_0 > 0$ such that,

$$(1/(2r^*))\rho(P^*(v_j)) \leq C_0(1/r_1)\rho(P(v_j)),$$

and

$$C_0(1/(2r^*))\rho(P^*(v_j)) \geq (1/r_{2k})\rho(P(v_j)) = \delta^{-(2k-1)}(1/r_1)\rho(P(v_j)),$$

where $j = 0, 1$. As $P^*(v_0) = P(v_0)$, it follows that $1/(2r^*) \leq C_0/r_1$ and $C_0/(2r^*) \geq \delta^{-(2k-1)}/r_1$. Then,

$$\rho(P^*(v_1))/\rho(P(v_1)) \leq C_0(2r^*)/r_1 \leq (C_0)^2\delta^{2k-1} = C,$$

and

$$\rho(P^*(v_1))/\rho(P(v_1)) \geq \delta^{-(2k-1)}(2r^*)/(C_0 r_1) \geq \delta^{-(2k-1)}/(C_0)^2 = 1/C.$$

Since the constant $C$ does not depend on $v_1$, Lemma 4.1 is proved in case $\tau(P) < \infty$. □

*Case* 2: $\tau(P) = \infty$. Let $W$ be the subset of vertices $v$ for which

$$\tau(P(v)) = \rho(P(v))/\text{dist}\,(0, P(v)) \geq 1.$$

Then clearly $W$ is infinite.

LEMMA 6.3. *Assume that $\tau(P) = +\infty$. Then*

$$\text{(6.13)} \qquad \lim_{v \in W, v \to \infty} \sup \tau(P^*(v)) > 0.$$

*Proof.* Assume the contrary. Then for any $0 < \varepsilon < 1/2$, there is a finite subset of vertices $V_0$ in $G$, with $v_0 \in V_0$, such that

$$\text{(6.14)} \qquad \rho(P^*(v))/\text{dist}\,(0, P^*(v)) = \tau(P^*(v)) < \varepsilon, \quad \text{for all } v \in W - V_0.$$

The choice of $\varepsilon$ will be made later. Let $W' = W - V_0$, and let $\tilde{G}$ be the subgraph of $G$ obtained by removing the vertices in $W'$ as well as the edges with an end in $W'$.

Let $r_0 > 0$ be such that $P(V_0) \subset D(r_0)$. Let $r_i = 4^i r_0$, and let $V_i = V_{c(r_i)}$, $i \geq 1$. Then for $i \neq j$, $V_i \cap V_j \subset W' = W - V_0$, and for any $0 \leq i_1 < i_2 < i_3$, $V_{i_2}$ separates $V_{i_1}$ from $V_{i_3}$ in $G$. Let $\text{VEL}_{\tilde{G}}(V_i, V_j)$ denote the vertex extremal length between $V_i \cap (V - W')$ and $V_j \cap (V - W')$ in the subgraph $\tilde{G}$.

Let $k$ be the following integer:

$$\text{(6.15)} \qquad k = 2 + \text{the integer part of } 2(128 + (16\pi)^2)(288 + (24\pi)^2).$$

As before, by applying Lemma 5.3 to $\tilde{G}$, we have,

$$(6.16) \qquad \mathrm{VEL}_{\tilde{G}}\left(V_2, V_{2k-1}\right) \geq \sum_{i=1}^{k-1} \mathrm{VEL}_{\tilde{G}}\left(V_{2i}, V_{2i+1}\right)$$

$$\geq \frac{k-1}{128 + (16\pi)^2} > 2[288 + (24\pi)^2].$$

We claim that there is some $r^* > 0$ such that for any $\rho \in [r^*, 2r^*]$, the vertex subset $V_{c(\rho)}^*$ separates $V_2$ from $V_{2k-1}$ in $G$. Assuming this, it follows that $V_2 \cap V_{2k-1} \subset V_{c(r^*)}^* \cap V_{c(2r^*)}^*$. On the other hand, $V_2 \cap V_{2k-1} \subset W' = W - V_0$. Then by (6.14), we deduce that $V_2 \cap V_{2k-1} \subset V_{c(r^*)}^* \cap V_{c(2r^*)}^* \cap (W - V_0) = \emptyset$. That is, $V_{c(4^2 r_0)} \cap V_{c(4^{2k-1} r_0)} = \emptyset$. As $r_0$ is arbitrary but sufficiently big, this implies that $\tau(P) \leq (4^{2k-3} - 1)/2$, a contradiction.

To show the existence of $r^*$ with the above property, let $U_{2,2k-1}$ be the set of vertices in $G$ which is separated from $V_0$ by $V_2$, and separated from $\infty$ by $V_{2k-1}$. In fact, $U_{2,2k-1}$ is equal to the set of vertices $v$ for which $P(v)$ has nonempty intersection with the annulus bounded by $c(r_2)$ and $c(r_{2k-1})$. Let $W'' = U_{2,2k-1} \cap W$. Then, clearly the cardinality of $W''$ is bounded by a universal constant $C_1$.

Let $r^* > 0$ be the minimum number such that $P^*(v) \cap D(r^*) \neq \emptyset$ for any $v \in V_2$. We may assume that $r^* = 1$. Let $\gamma^* \in \Gamma_G^*(V_2, V_{2k-1})$. As in Lemma 6.1, the diameter of $P^*(\gamma^*)$ is at least equal to $r^*$:

$$(6.17) \qquad \mathrm{diam}\left(P^*(\gamma^*)\right) \geq r^* = 1.$$

By contradiction, let us assume that there is a circle $c(\rho_1)$ with $\rho_1 \in [1, 2]$ such that $V_{c(\rho_1)}^*$ does not separate $V_2$ from $V_{2k-1}$. Thus there is a path $\gamma_0$ in the graph joining $V_2$ and $V_{2k-1}$ such that $\gamma_0 \cap V_{c(\rho_1)}^* = \emptyset$. As in Lemma 6.1, $P^*(\gamma_0)$ is contained in the interior of $D(\rho_1)$, and then,

$$(6.18) \qquad P^*(\gamma^*) \cap D(\rho_1) \supseteq P^*(\gamma^* \cap \gamma_0) \neq \emptyset.$$

For any vertex $v \in V$, let $\eta(v) = \mathrm{diam}\left(P^*(v) \cap D(3)\right)$. Using (6.17) and (6.18), and since $3 \geq \rho_1 + 1$, we deduce that for any *connected* vertex curve $\gamma^*$ in $\Gamma_G^*(V_2, V_{2k-1})$,

$$(6.19) \qquad \sum_{v \in \gamma^*} \eta(v) \geq 1.$$

But since $V_2$ and $V_{2k-1}$ are finite and connected and $G$ is a disk triangulation graph, any vertex curve in $\Gamma_G^*(V_2, V_{2k-1})$ contains a *connected* vertex subcurve in $\Gamma_G^*(V_2, V_{2k-1})$. Thus (6.19) holds for any $\gamma^*$ in $\Gamma_G^*(V_2, V_{2k-1})$.

By (6.14), we have $\eta(v) \leq 6\varepsilon$ for any $v \in W' \supseteq W''$. Since the cardinality of $W''$ is at most $C_1$, we see that

$$(6.20) \qquad \sum_{v \in W''} \eta(v) \leq 6C_1 \varepsilon.$$

Now let $\beta^* \in \Gamma^*_{\tilde{G}}\big(V_2 \cap (V - W'), V_{2k-1} \cap (V - W')\big)$. Then, the union $\gamma^* = \beta^* \cup W''$ is a vertex curve in $\Gamma^*_G(V_2, V_{2k-1})$. Using (6.19) and (6.20), we deduce that

$$\sum_{v \in \beta^*} \eta(v) \geq 1 - 6C_1\varepsilon.$$

We choose $\varepsilon > 0$ so small that $1 - 6C_1\varepsilon \geq 1/\sqrt{2}$. Then $\sqrt{2}\eta$ is $\Gamma^*_{\tilde{G}}\big(V_2 \cap (V - W'), V_{2k-1} \cap (V - W')\big)$-admissible. Thus by (6.1),

$$\begin{aligned}
\mathrm{VEL}_{\tilde{G}}(V_2, V_{2k-1}) &= \mathrm{MOD}\left(\Gamma^*_{\tilde{G}}(V_2 \cap (V - W'), V_{2k-1} \cap (V - W'))\right) \\
&\leq \mathrm{area}\,(\sqrt{2}\eta).
\end{aligned}$$

As in Lemma 6.1, it follows that

$$\mathrm{VEL}_{\tilde{G}}(V_2, V_{2k-1}) \leq 2 \cdot \mathrm{area}\,(\eta) \leq 2[32 + (8\pi)^2] \cdot 3^2.$$

This contradicts (6.16). The proof of Lemma 6.3 is thus complete. $\qquad\square$

*Proof of Lemma* 6.1, *Case* 2 (*continued*). Let $\delta \in (0, 1/3)$ be a fixed number strictly smaller than third of the limit in (6.13). Then, there is a sequence of mutually distinct vertices $u_k$ such that $\tau(P(u_k)) \geq 1$ and $\tau(P^*(u_k)) \geq 3\delta$. As $P$ and $P^*$ are both locally finite in $\mathbb{C}$ (see Cor. 6.2), we have $\mathrm{dist}\,(0, P(u_k)) \to \infty$ and $\mathrm{dist}\,(0, P^*(u_k)) \to \infty$.

Let $\rho_0 > 0$ be such that $P(v_0) \cup P(v_1) \cup P^*(v_0) \cup P^*(v_1) \subseteq D(\rho_0)$. Let $u = u_k$ be a vertex in the above sequence such that $d = \mathrm{dist}\,(0, P(u_k)) \geq 100\rho_0/\delta^2$ and $d^* = \mathrm{dist}\,(0, P^*(u_k)) \geq 100\rho_0/\delta^2$. See Figure 6.3 for the configuration of $D(\rho_0)$ and $P(u)$. The figure of $P^*(u)$ is similar.

Let $F$ be a Möbius transformation with $F(0) = 0$ and $F(\hat{\mathbb{C}} - P(u)) = U$, and let $F^*$ be a Möbius transformation with $F^*(0) = 0$ and $F^*(\hat{\mathbb{C}} - P^*(u)) = U$. Then the hyperbolic distance between $F(\infty)$ and $0$ is bigger than the hyperbolic distance between $0$ and $\rho(P(u))/(d + \rho_0 + \rho(P(u)))$. It follows then $|F(\infty)| = \rho(P(u))/(d + \rho(P(u))) = \tau(P(u))/(1 + \tau(P(u))) > 1/3 > \delta$. Similarly, $|F^*(\infty)| = \tau(P^*(u))/(1 + \tau(P^*(u))) \geq 3\delta/(1.01 + 3\delta) > \delta$. Also, it is easy to see from $\rho_0 \leq \delta^2 d/100$ and $\rho_0 \leq \delta^2 d^*/100$ that the disks $F(D(\rho_0))$, $F^*(D(\rho_0))$; hence $F(P(v_0))$, $F(P(v_1))$, $F^*(P^*(v_0))$, and $F^*(P(v_1))$ are all contained in $(\delta/2)U$.

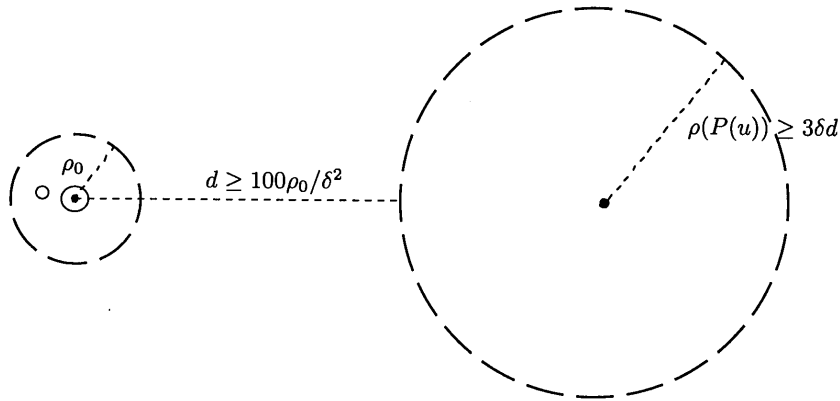

**Figure 6.3.** The disks $D(\rho_0)$ and $P(u)$ (not drawn to scale).

Since dist $(0, P(u)) \geq 100\rho_0$, the absolute derivative $|dF/dz|$ of $F$ in $D(\rho_0)$ does not vary by a factor $\geq 2$. Therefore, the ratio between

$$\rho(F(P(P(v_1))))/\rho(F(P(P(v_0))))$$

and

$$\rho(P(v_1))/\rho(P(v_0))$$

is in $[1/4, 4]$. Similarly, the ratio between

$$\rho(F^*(P^*(v_1)))/\rho(F^*(P^*(v_0)))$$

and

$$\rho(P^*(v_1))/\rho(P^*(v_0))$$

is in $[1/4, 4]$. As $P(v_0) = P^*(v_0)$, in order to get lower and upper bounds for $\rho(P^*(v_1))/\rho(P(v_1))$, we need show that $\rho(F(P(v_1)))/\rho(F(P(v_0)))$ is comparable with $\rho(F^*(P^*(v_1)))/\rho(F^*(P^*(v_0)))$.

Compare the pattern $(1/\delta)F(P)$ with $F^*(P^*)$. In the first pattern, as $|F(\infty)| > \delta$ and $F(P)$ is locally finite in $\hat{\mathbb{C}} - \{F(\infty)\}$, only a finite number of disks intersect with $U$. Using Lemma 2.3 as in Case 1, we obtain, for $j = 0, 1$,

$$\rho_{\mathrm{hyp}}\left((1/\delta)F(P(v_j))\right) \geq \rho_{\mathrm{hyp}}\left(F^*(P^*(v_j))\right).$$

Similarly, comparing $F(P)$ with $(1/\delta)F^*(P^*)$, as $|F^*(\infty)| > \delta$, we obtain

$$\rho_{\mathrm{hyp}}\left(F(P(v_j))\right) \leq \rho_{\mathrm{hyp}}\left((1/\delta)F^*(P^*(v_j))\right).$$

Lemma 4.1 follows since all the disks involved are contained in $(1/2)U$.    $\square$

## 7. Uniformization of disk patterns

In this section, we will prove Uniformization Theorem 1.3. The following is an analog of the Ring Lemma of [21]:

LEMMA 7.1. *Let $G$ be a graph and let $\Theta : E \to [0, \pi/2]$ be a function on the set of edges. Let $P$ be a disk pattern in $\mathbb{C}$ which realizes the data $(G, \Theta)$. Let $v_0$ be an interior vertex, and let $v_1, v_2, \ldots, v_l, v_{l+1} = v_1$ be a closed chain of neighboring vertices. Assume that the vertices $v_k$, $1 \le k \le l$, are also interior vertices. Then there is a constant $C_1 = C_1(v_0, G, \Theta)$ depending on $v_0$, $G$, and $\Theta$, such that,*

$$(7.1) \qquad \rho(P(v_0)) \le C_1 \rho(P(v_k)), \quad 1 \le k \le l.$$

*Proof.* Assume that no such constant $C_1$ exists. Then there is a sequence of patterns $P_n$ in the plane which realize $(G, \Theta)$, such that $P_n(v_0) = D(1)$, and $\lim_{n \to \infty} \rho(P_n(v_k)) = 0$ for some $k$, $1 \le k \le l$. By subtracting a subsequence, we may assume that for each vertex $v$ in the graph, $\lim_{n \to \infty} \rho(P_n(v))$ exists in $[0, \infty]$. Denote by $W_+$ the set of vertices $v$ in the graph for which $\lim_{n \to \infty} \rho(P_n(v)) > 0$, and $W_- = V - W_+$.

For any integer $m \ge 0$, let $U_m$ be the set of vertices in $G$ whose combinatorial distance from $v_0$ is at most $m$, and let $V_m = U_m - U_{m-1}$. Then clearly, $V_1$ contains at least two vertices in $W_+$ and $V_2$ contains at least three vertices in $W_+$. We have $v_k \in V_-$. Let $W'_-$ be the set of vertices in $U_2 \cap W_-$, which may be joined to $v_k$ by a path passing through vertices in $U_2 \cap W_-$ only. Let $W'_+$ be the set of vertices in $U_2$ which share an edge with a vertex in $W'_-$. It is then easy to see that $W'_+$ consists of the (distinct) vertices: $v_0, v_{1,i}, 1 \le i \le i_1$, and $v_{2,j}, 1 \le j \le j_1$, where $v_{1,i} \in V_1$, $v_{2,j} \in V_2$, $i_1 \ge 2$ and $j_1 \ge 1$.

Replacing by a subsequence if necessary, we may assume that for each $v \in W'_+$, the sequence of disks $P_n(v)$ converges to some disk $P_\infty(v)$ in $\hat{\mathbb{C}}$. These limit disks all pass through a single point on the unit circle $\partial P_\infty(v_0) = c(1)$. Moreover, for each $j$, as the vertices $v_0$ and $v_{2,j}$ do not share an edge in $G$, the disks $P_\infty(v_0)$ and $P_\infty(v_{2,j})$ should then be tangent. It follows that $j_1 = 1$; otherwise $\Theta(v_{2,1}v_{2,2})$ would be equal to $\pi$. Similarly, $i_1 = 2$. Then $w_0 = v_0$, $w_1 = v_{2,1}$, $w_2 = v_{2,1}$, $w_3 = v_{1,2}$, $w_4 = w_0$ form a simple loop of four edges in $G$ which separates the vertices of $G$. Moreover, $\Theta(w_i w_{i+1}) = \pi/2$. This contradicts conditions (C1) and (C2) in Section 1. The lemma is thus proved. $\square$

Let $G$ be a disk triangulation graph and let $\Theta : E \to [0, \pi/2]$ be a function on the set of edges. Suppose that conditions (C1) and (C2) in Section 1 are satisfied. Let $G_n$ be an increasing sequence of subgraphs of $G$ whose union is $G$. Again, we will choose $G_n$ to be one-skeletons of triangulations of the closed topological disk. Then using Andreev's theorem ([26]), for each $n$, there is a disk pattern $P_n$ which realizes $(G_n, \Theta|_{G_n})$. Moreover, we may assume that all disks are contained in the closed unit disk $\overline{U} = D(1)$, and that the boundary disks are all internally tangent to $\partial U$.

Let $v_0$ be a fixed vertex. By means of Möbius transformations preserving $\overline{U}$, we may assume that for each $n$, the disk $P_n(v_0)$ is centered at 0. Let $\delta_n = 1/\rho(P_n(v_0))$. Then by Lemma 2.3, $\delta_n$ is an increasing sequence. Consider $Q_n = \delta_n P_n$. Then $Q_n(v_0) = D(1)$ for each $n$. Let $v$ be a vertex in $G$. Using Lemma 7.1 inductively, there is a constant $C_2(v) = C_2(v, G, \Theta) \geq 1$ such that for any big $n$,

$$(7.2) \qquad\qquad 1/C_2(v) \leq \rho(Q_n(v)) \leq C_2(v).$$

Thus by subtracting a subsequence if necessary, we may assume that each sequence $Q_n(v)$ converges to a bounded disk $Q(v)$ in $\mathbb{C}$. Clearly, the dihedral angle of a pair of neighboring disks is preserved in the limit. Thus we have proved the following lemma:

LEMMA 7.2.    *Let $G$ be a disk triangulation graph and let $\Theta : E \to [0, \pi/2]$ be a function of edges. Suppose that conditions* (C1) *and* (C2) *are satisfied. Then there is a disk pattern $Q$ in $\mathbb{C}$ which realizes the data $(G, \Theta)$.*

*Remark.* With a similar proof, Lemma 7.2 also holds if $G$ is the one-skeleton of a triangulation of a planar surface possibly with boundary so that each boundary component has at least four vertices, or just any connected planar graph with at least five vertices such that each edge has two distinct vertices and through a pair of vertices there is at most one edge.

Theorem 1.3 follows from Lemma 7.2, Corollary 6.2, and the following lemma.

LEMMA 7.3.    *Let $G$ be a disk triangulation graph and let $\Theta : E \to [0, \pi/2]$. Suppose that conditions* (C1) *and* (C2) *are satisfied. Then the following statements are equivalent*:

(1) *The graph $G$ is* VEL-*hyperbolic*;
(2) *For any increasing sequence of finite, connected subsets of vertices $V_k$ in $G$ such that $V_\infty = \bigcup_{k=1}^{\infty} V_k$ is infinite, we have $\lim_{k \to \infty} \mathrm{VEL}\,(V_k, \infty) = 0$;*
(3) *There is a disk pattern $P$, locally finite in $U$, which realizes $(G, \Theta)$.*

*Proof.* By Corollary 6.2, (3) implies (1). Assuming that (1) holds, let us prove (2). As the sequence $\mathrm{VEL}\,(V_k, \infty)$ is monotone decreasing, we need only show that a subsequence converges to 0. Let $Q$ be a disk pattern in $\mathbb{C}$ which realizes $(G, \Theta)$ (actually we may assume $\Theta = 0$). There is a unique simply connected domain $\Omega$ in $\mathbb{C}$, such that $Q$ is (contained and) locally finite in $\Omega$. By Corollary 6.2, $\Omega \neq \mathbb{C}$. Let $x_\infty \in \partial\Omega$ be a limit point of $Q(V_\infty) =$ $\bigcup_{v \in V_\infty} Q(v) \subseteq \Omega$. Inverting the disk pattern on a circle centered at $x_\infty$, we may assume that $x_\infty = \infty$. Let $r_0 > 0$ be big such that $Q(V_1) \subset D(r_0)$ and that $D(r_0) \cap \partial\Omega \neq \emptyset$. Let $r_k$, $k \geq 1$, be an increasing sequence such that $r_k \geq 2r_{k-1}$, and that $V_{c(r_k)} \cap V_{c(r_{k-1})} = \emptyset$, where $V_c = V_c(Q) = \{v; \ Q(v) \cap c \neq \emptyset\}$.

Let $s_k \geq 0$ be the minimum number such that $D(s_k) \cap Q(v) \neq \emptyset$ for any $v \in V_k$. Since $\infty$ is a limit point of $Q(V_\infty)$, we deduce that $s_k$ diverges to $\infty$. Dropping a subsequence of $V_k$ if necessary, we may assume that $s_k \geq r_k$. Let $k = 2m$. Then, as $D(r_0) \cap \partial\Omega \neq \emptyset$, any vertex curve in $\Gamma_G^*(V_{2m}, \infty)$ contains a path connecting $V_{c(r_1)}$ to $V_{c(r_{2m})}$, and hence a disjoint union of subpaths $\gamma_i$, $1 \leq i \leq m$, with $\gamma_i \in \Gamma_G(V_{c(r_{2i-1})}, V_{c(r_{2i})})$. It follows that

$$\mathrm{VEL}\left(\Gamma_G^*(V_{2m}, \infty)\right) \geq \sum_{i=1}^m \mathrm{VEL}\left(V_{c(r_{2i-1})}, V_{c(r_{2i})}\right).$$

Then by Lemma 5.3, we deduce that,

$$\mathrm{VEL}\left(\Gamma_G^*(V_{2m}, \infty\right) \geq \frac{m}{128 + (16\pi)^2}.$$

Thus by (6.1), we obtain $\mathrm{VEL}(V_{2m}, \infty)) = 1/\mathrm{VEL}\left(\Gamma_G^*(V_{2m}, \infty)\right) \to 0$.

It remains to show that (2) implies (3). As $\mathrm{VEL}(V_k, \infty)$ is finite for sufficiently big $k$, the graph $G$ is VEL-hyperbolic. Let $G_n$, $P_n$, $\delta_n = 1/\rho(P_n(v_0))$, $Q_n = \delta_n P_n$ and $Q = \lim_{n \to \infty} Q_n$ be as constructed before the statement of Lemma 7.2. Let $\delta_\infty = \lim_{n \to \infty} \delta_n$. Clearly, $Q$ is contained in $D(\delta_\infty)$, and hence in $\mathrm{int}\, D(\delta_\infty) = \{z \in \mathbb{C}; \ |z| < \delta_\infty\}$ because $G$ is a disk triangulation graph. We will show that $Q$ is locally finite in $\mathrm{int}\, D(\delta_\infty)$. Then it follows by Corollary 6.2 that $\delta_\infty$ is finite and therefore $P = (1/\delta_\infty)Q$ is the required disk pattern.

By contradiction, assume that $Q$ is not locally finite in $\mathrm{int}\, D(\delta_\infty)$. Let $\delta_0' < \delta_\infty$ be such that the set $W$ of vertices $v$ for which $Q(v) \cap D(\delta_0') \neq \emptyset$ is infinite. Then it is elementary to see that $W$ contains an infinite simple path $(v_0, v_1, v_2, \dots, v_k, \dots)$ in the graph. Let $V_k = \{v_0, v_1, \dots, v_k\}$. Let $\delta_1'$, $\delta_2'$, be some fixed numbers with $\delta_0' < \delta_1' < \delta_2' < \delta_\infty$. Then for each $k$, there is some $n = n_k$ such that $\delta_{n_k} > \delta_2'$ and $Q_{n_k}(v_i) \cap D(\delta_1') \neq \emptyset$ for $0 \leq i \leq k$. Applying Lemma 5.3 to $Q_{n_k}$, we obtain,

$$\mathrm{VEL}\left(V_k, \partial G_{n_k}\right) \geq \mathrm{VEL}\left(V_{c(\delta_1')}(Q_{n_k}), V_{c(\delta_2')}(Q_{n_k})\right) \geq \frac{(\delta_2' - \delta_1')^2}{(\delta_2')^2(32 + (8\pi)^2)},$$

where $\partial G_{n_k}$ denotes the set of boundary vertices of $G_{n_k}$. In particular, $\mathrm{VEL}(V_k, \infty) \geq \mathrm{VEL}(V_k, \partial G_{n_k})$ cannot converge to 0, a contradiction to (2). $\square$

*Remark.* We do not know if the equivalence of (1) and (2) in the above lemma holds for any locally finite, connected graph with a finite number of ends.

UNIVERSITY OF CALIFORNIA, SAN DIEGO, LA JOLLA, CA
*E-mail address*: zhe@math.ucsd.edu